\documentclass[11pt,a4paper,reqno]{amsart}

\usepackage[T1]{fontenc}
\usepackage{microtype, a4wide, textcomp,fbb}
\usepackage{
    amsmath,  amssymb,  amsthm,   amscd,
    gensymb,  graphicx, etoolbox, 
    booktabs, stackrel, mathtools    
}
\usepackage[usenames,dvipsnames]{xcolor}
\definecolor{darkblue}{rgb}{0,0,0.7}
\definecolor{darkred}{rgb}{0.7,0,0}
\usepackage[colorlinks=true, linkcolor=darkred, citecolor=darkblue, urlcolor=blue, pagebackref=true, breaklinks=true]{hyperref}
\usepackage[capitalise]{cleveref}
\usepackage{enumitem}
\usepackage{placeins}
\usepackage{relsize}
\usepackage{todonotes}
\usepackage{soul}
\usepackage{tikz}
\usetikzlibrary{arrows}
\usetikzlibrary{shapes}
\usepackage{float}
\usepackage{tabularx}
\usepackage{kantlipsum}
\usepackage{array}
\usepackage[margin=2.6cm]{geometry}

\newtheorem{proposition}{Proposition}[section]
\newtheorem{lemma}[proposition]{Lemma}
\newtheorem{theorem}[proposition]{Theorem}
\newtheorem{corollary}[proposition]{Corollary}
\newtheorem{question}[proposition]{Question}
\newtheorem{conjecture}[proposition]{Conjecture}
\newtheorem{remark}[proposition]{Remark}

\newtheorem{definition}[proposition]{Definition}

\newtheorem{innercustomthm}{Theorem}
\newenvironment{customthm}[1]
  {\renewcommand\theinnercustomthm{#1}\innercustomthm\itshape}
  {\endinnercustomthm}

\newcommand{\reg}{{\rm reg}}

\def\H{\mathcal{H}}

\tikzstyle{place}=[draw,circle,minimum size=1mm,inner sep=1pt,outer sep=-1.1pt,fill=black]

\usetikzlibrary{shapes.geometric}
\tikzstyle{places}=[draw,rectangle,minimum size=8pt,inner sep=0pt]
\tikzstyle{placesf}=[draw,rectangle,minimum size=5pt,inner sep=0pt]
\tikzstyle{placec}=[draw,circle,minimum size=8pt,inner sep=0pt]
\tikzstyle{placecf}=[draw,circle, minimum size=5pt,inner sep=0pt]

\def\K{\mathbb{K}}
\def\reg{\mathrm{reg}}
\def\lcm{\mathrm{lcm}}
\def\level{\mathrm{level}}

\def\pd{\mathrm{pd}}

\def\H{\mathcal H}
\def\P{\mathcal P}
\def\E{\mathcal E}
\def\F{\mathcal F}
\def\G{\mathcal G}
\def\x{\mathbf x}

\def\l{\langle}
\def\r{\rangle}
\def\b{\mathrm{bight}}
\def\d{\mathrm{depth}}
\def\p{\mathfrak p}
\def\q{\mathfrak q}
\def\h{\mathrm{ht}}

\begin{document}

\title{On the path ideals of chordal graphs}

\author{Kanoy Kumar Das}
\address{Chennai Mathematical Institute, India}
\email{kanoydas@cmi.ac.in; kanoydas0296@gmail.com}

\author{Amit Roy}
\address{Chennai Mathematical Institute, India}
\email{amitiisermohali493@gmail.com}

\author{Kamalesh Saha}
\address{Chennai Mathematical Institute, India}
\email{ksaha@cmi.ac.in; kamalesh.saha44@gmail.com}

\keywords{$t$-path ideal, chordal graph, tree, Alexander dual, regularity, projective dimension}
\subjclass[2020]{05E40, 13F55,  05C05}

\vspace*{-0.4cm}
\begin{abstract}
%The regularity and projective dimension of edge ideals of chordal graphs are well-known in terms of the induced matching number and the big height, respectively. In this paper, we extend these results for $3$-path ideals of chordal graphs. As a consequence, we get that the $3$-path ideal of a chordal graph is Cohen-Macaulay if and only if it is unmixed. Further, we show the Alexander dual of the $3$-path ideal of a tree is vertex splittable, and this does not hold for all chordal graphs. This settles the $t=3$ case of a recent conjecture in \cite{ADGRS}. Furthermore, we extend the formula of the regularity of the $3$-path ideal of chordal graphs for all $t$-path ideals of caterpillar graphs. We then provide some families of graphs to show that these formulas for the regularity and the projective dimension cannot be extended to higher $t$-path ideals of chordal graphs, or even for trees. 

In this article, we give combinatorial formulas for the regularity and the projective dimension of $3$-path ideals of chordal graphs, extending the well-known formulas for the edge ideals of chordal graphs given in terms of the induced matching number and the big height, respectively. As a consequence, we get that the $3$-path ideal of a chordal graph is Cohen-Macaulay if and only if it is unmixed. Additionally, we show that the Alexander dual of the $3$-path ideal of a tree is vertex splittable, thereby resolving the $t=3$ case of a recent conjecture in \cite{ADGRS}. Also, we give examples of chordal graphs where the duals of their $t$-path ideals are not vertex splittable for $t\ge 3$. Furthermore, we extend the formula of the regularity of $3$-path ideals of chordal graphs to all $t$-path ideals of caterpillar graphs. We then provide some families of graphs to show that these formulas for the regularity and the projective dimension cannot be extended to higher $t$-path ideals of chordal graphs (even in the case of trees).

%[Internat. J. Algebra Comput., 33(3):481–498, 2023]

 % We show that for chordal graphs, the Castelnuovo-Mumford regularity and the projective dimension of the 3-path ideal are the same as \textcolor{blue}{the $3$-path induced matching number of the graph $G$,} and the big height of the $3$-path ideal, respectively. As a consequence, we get that the $3$-path ideal of a chordal graph is Cohen-Macaulay if and only if it is unmixed.  Furthermore, we show that the formula for the regularity of the $3$-path ideal can be extended to the $t$-path ideals of caterpillar graphs. We then provide a class of examples of $t$-path ideals to show that the above formulas for the regularity and the projective dimension cannot be extended to higher $t$-path ideals of chordal graphs, or even for trees. We finally study the Alexander dual of the $3$-path ideals of trees and show that these are vertex splittable. This settles the $t=3$ case of a recent conjecture in \cite{ADGRS}. 
\end{abstract}

\maketitle

\section{Introduction}

Let $G$ be a simple graph with the set of vertices $V(G)$ and the set of edges $E(G)$. Identifying the vertices of $G$ as the indeterminates, we consider the polynomial ring $R=\K[v:v\in V(G)]$, where $\K$ is any field. There are several ways to associate the graph $G$ with an ideal in $R$. The first and perhaps the most famous of them is given by Villarreal \cite{CM} in 1990, which is known as the edge ideal of $G$, denoted by $I(G)$. Edge ideals of graphs are quadratic square-free monomial ideals and a significant portion of the literature is dedicated to the exploration of such ideals.
The primary objective of these studies is to establish a connection between the combinatorial structures of a graph and the algebraic properties of the corresponding ideal. The Castelnuovo-Mumford regularity (in short, regularity) is an algebraic invariant which somewhat measures the complexity of a graded ideal or a graded module and can be regarded as a crude estimation of the Betti numbers. Another such algebraic invariant is the projective dimension, which measures the length of a minimal free resolution. Interestingly, for various well-known classes of graphs, the regularity and the projective dimension of the edge ideal can be estimated by combinatorial invariants of the associated graphs (See for instance \cite{Ha2014, MoreyVillarreal}). Although, finding a unified formula for these invariants remains an open question.  
This article focuses on exploring these two algebraic invariants for a special class of monomial ideals, namely the $t$-path ideals.

\par 

A generalization of edge ideals is the so-called t-path ideals, first introduced by Conca, and De Negri \cite{ConcaDeNegri1998} for directed graphs. Note that for $t=2$,  
they are the usual edge ideals of graphs. Later,  Alilooee and Faridi \cite{AlilooeeFaridi2015, AlilooeeFaridi2018} considered $t$-path ideals of simple undirected graphs, which we shall denote as $I_t(G)$. They studied Betti numbers of $t$-path ideals of lines and cycles and gave precise formulas for the regularity and the projective dimension. In this article, we are interested in finding the regularity and projective dimension of the ideals $I_t(G)$, when $G$ is a chordal graph. We now briefly recall the results which are already known in the literature related to regularity and projective dimension of $t$-path ideals. Banerjee \cite{Banerjee2017} showed that if $G$ is gap-free and claw-free, then $I_t(G)$ has a linear resolution for $t=3,4,5,6$. Moreover, it was shown in \cite{Banerjee2017} that $I_t(G)$ has a linear resolution for all $t\ge 3$ if $G$ is gap-free, claw-free and whiskered $K_4$-free. In \cite{Erey2020}, Erey gave a formula for the Betti numbers of path ideals of line and star graphs. Also, Bouchat, Ha, and O’Keefe \cite{BouchatHaOkeefe2011}, He and Van Tuyl \cite{HeVanTuyl2010}, Kiani and Madani \cite{KianiMadani2016}, and Bouchat and Brown \cite{BouchatBrown2017} studied path ideals of directed rooted trees, which are in general subideals of the path ideals of the underlying undirected graphs.
\par

In a recent article \cite{DRSV23}, Deshpande, Singh, and Van Tuyl, along with the second author, studied another generalization of edge ideals, which is similar to the t-path ideals and is known as the t-connected ideals, denoted by $J_t(G)$ (see also \cite{AJM2024}). These ideals can be interpreted as Stanley-Reisner ideals of simplicial complexes constructed from the graph $G$. They used this description to give a partial generalization of Froberg's theorem in the case of $t$-connected ideals. More precisely, it has been proved in \cite{DRSV23} that if $G$ is co-chordal, then $J_t(G)$ is vertex splittable for all $t \geq 2$. It is worthwhile to mention that, for $t=3$, the 3-connected ideals coincide with the 3-path ideals, that is, $J_3(G)=I_3(G)$ for any simple graph $G$. Thus, the class of $3$-path ideals provides a shared platform to explore which properties of edge ideals of simple graphs can be generalized to higher $t$-path ideals or $t$-connected ideals. 

In the literature, there are very few results known for $3$-path ideals, when compared to edge ideals of graphs. A recent study conducted by Kumar and Sarkar \cite{KumarSarkar2024} shows that if $G$ is a tree, then $\reg(R/I_3(G)) = 2\nu_3(G)$, where $\nu_3(G)$ is the 3-path induced matching number of $G$ (see Section \ref{induced_matching} for the definitions). Later, Hang and Vu \cite{HangVu2024} gave a complete formula for the $\reg(R/I_3(G))$ and $\pd(R/I_3(G))$, when $G$ is a connected unicyclic graph. To the best of our knowledge, except when $I_3(G)$ has a linear resolution, an exact formula for the regularity or projective dimension of the $3$-path ideals is known only if $G$ is a tree or a unicyclic graph.

In this work, we show that the regularity formula of $3$-path ideals for trees in \cite{KumarSarkar2024} also holds for the class of chordal graphs. In addition, we also give a formula for the projective dimension of the $3$-path ideals of chordal graphs. However, our approach is slightly different from the previous known works related to path ideals. A popular technique in this context is to make use of inductive methods by splitting a monomial ideal $I$ into two smaller ideals $J$ and $K$ such that $I=J+K$, where the minimal generating sets of $J$ and $K$ are disjoint. We make use of a similar splitting for $I_3(G)$ by means of a particular edge in the underlying chordal graph $G$, omitting the condition on the minimal generating sets of the subideals. Using this new splitting of $I_3(G)$, we prove the following:

\vskip 0.2cm
\noindent 
\textbf{Theorem} (\cref{result1,pdmain}).
Let $G$ be a chordal graph. Then 
\begin{enumerate}
    \item $\reg(R/I_3(G))=2\nu_3(G)$, where $\nu_3(G)$ denotes the $3$-path induced matching number of $G$.

    \item $\pd(R/I_3(G))=\b(I_3(G))$,    where $\b(I_3(G))$ denotes the maximum height of a minimal prime ideal of $I_3(G)$.
\end{enumerate}

As a consequence of the above result, we obtain an equivalent condition for the Cohen-Macaulayness of the $3$-path ideals of chordal graphs in \cref{CM}, extending the works of \cite{CamposGundersonMoreyPaulsenPolstra2014} for the case $t=3$. More precisely, we show that for a chordal graph $G$, $I_3(G)$ is Cohen-Macaulay if and only if it is unmixed.

Our second main result uses some inspiration from the works of \cite{DRSV23}. In particular, we show the following:

\begin{customthm}{\ref{main theorem1}}
    Let $T$ be a tree.
    Then the Alexander dual ideal $I_3(T)^{\vee}$ is a vertex splittable ideal.
\end{customthm}

The above result has several implications, and one of them is that the ideal $I_3(T)$ is sequentially Cohen-Macaulay, and this property does not depend on the characteristic of the base field $\K$. Furthermore, a direct consequence of \Cref{main theorem1} is that the Stanley-Reisner complex of the ideal $I_3(T)$ is vertex decomposable, which settles \cite[Conjecture 3.15]{ADGRS} for $t=3$ case (see \Cref{vdtree}). \par 

One of the most crucial parts of this article is \Cref{section 6}, where we discuss how far one should investigate in generalizing the results related to edge ideals and $3$-path ideals for higher $t$-path ideals. For instance, is it possible to extend the formulas for the regularity and projective dimension given in \Cref{result1} and \Cref{pdmain}, respectively, for all $t$-path ideals of chordal graphs? In this regard, we show that for the class of caterpillar graphs, a similar regularity formula holds for all $t$-path ideals. More precisely,

\begin{customthm}{\ref{caterpillar reg}}
    Let $G$ be a caterpillar graph. Then for all $t\ge 2$, $\reg(R/I_t(G))= (t-1)\nu_t(G)$, where $\nu_t(G)$ denotes the $t$-path induced matching of $G$.
\end{customthm}

\noindent
However, for a general tree, we show in \Cref{reg tree} that one can not expect such a formula of the regularity for higher $t$-path ideals when compared to \Cref{result1} and \Cref{caterpillar reg}. Specifically,

\begin{customthm}{\ref{reg tree}}
    For a tree $T$ and for each $t\ge 4$, the difference between $\reg(R/I_t(T))$ and the quantity $(t-1)\nu_t(T)$ can be arbitrarily large for each $t\geq 4$.
\end{customthm}

\noindent
The above result disproves a recent conjecture \cite[Conjecture 4.9]{HangVu2024} by Hang and Vu for $t\ge 4$. Moving on, we show in \cref{Example1} that the vertex splittable property of the dual of $3$-path ideals of trees cannot be extended to the dual of $3$-path ideals of chordal graphs. 
%In fact, Also, the analogous formula for the projective dimension of $t$-path ideal does not hold even for the caterpillar graphs. 
In fact, we show in \cref{caterpillar not vs} that for a caterpillar graph $G$, the difference between $\pd(R/I_t(G))$ and $\b(I_t(G))$ can be arbitrarily large for $t\geq 4$. That is, in case of trees, the vertex splittable property of the dual of $3$-path ideals cannot be extended for the dual of higher $t$-path ideals.

%Thus, the algebraic invariant $\b(I_t(G))$ might not be a suitable parameter to estimate the projective dimension of $R/I_t(G)$ for $t\geq 4$. 

Any square-free monomial ideal can be viewed as an edge ideal of a hypergraph, as well as a Stanley-Reisner ideal of a simplicial complex. In this article, we make use of both the description whenever it is necessary.
Our paper is organised as follows. In \Cref{section 2}, we make precise definitions and fix some notations that will be used throughout the article. In this section we also list some of the known results that will be needed in the later sections. In \Cref{section 3}, we prove the formula for the regularity and the projective dimension of the 3-path ideals of chordal graphs. In \Cref{section 4}, we show that the Alexander dual of the 3-path ideal of a tree is vertex splittable. As a consequence of this result, we get that the 3-path ideal of a tree is sequentially Cohen-Macaulay. Finally, in \Cref{section 6}, we discuss some possible generalisations of all of our main results and propose some conjectures and questions related to higher $t$-path ideals and $t$-connected ideals. 

\section{Preliminaries}\label{section 2}

In this section, we recall some basic notions from combinatorics and commutative algebra, which will be used throughout the remaining part of the paper.

\subsection{Graph Theory and Combinatorics:}

    A (simple) {\it graph} $G$ is a pair $(V (G), E(G))$, where $V(G)$ is called the vertex set of $G$ and $E(G)$, a collection of cardinality $2$ subsets of $V(G)$, is known as the edge set of $G$. We now recall some useful notations that will be needed in the later sections.
    \begin{enumerate}
        \item If $x_1,\ldots,x_r\in V(G)$, then $G\setminus \{x_1,\ldots,x_r\}$ denotes the graph with the vertex set $V(G)\setminus \{x_1,\ldots,x_r\}$ and the edge set $\{\{u,v\}\in E(G)\mid x_i\notin \{u,v\}\text{ for each }i\in[r]\}$. If $r=1$, then $G\setminus\{x_1\}$ is simply denoted by $G\setminus x_1$.

        \item If $e_1,\ldots , e_s\in E(G)$, then $G- \{e_1,\ldots,e_s\}$ denotes the graph with the vertex set $V(G)$ and the edge set $\{\{u,v\}\in E(G)\mid e_i\neq \{u,v\}\text{ for each }i\in[s]\}$. As before, if $s=1$, then $G -\{e_1\}$ is simply denoted by $G - e_1$.

        \item For a vertex $a$ in $G$, the set of {\it neighbors} of $a$, denoted by $N_G(a)$, is the set $\{b \in V(G) \mid \{a,b\}\in E(G)\}$. 

        \item The set of \textit{closed neighbors} of $a$ is the set $N_G(a)\cup \{a\}$, and is denoted by $N_G[a]$. More generally, if $a_1,\ldots,a_k\in V(G)$ then we define $N_G[a_1,\ldots,a_k]=\{a_1,\ldots,a_k\}\cup\{b\in V(G)\mid \{b,a_i\}\in E(G)\text{ for some }1\leq i\leq k\}$.

        \item The number $|N_G(a)|$ is called the \textit{degree} of a vertex $a$ in $G$, and is denoted by $\deg(a)$. A {\it leaf} is a vertex of $G$ with degree exactly $1$.

        \item Let $W\subseteq V(G)$. Then the {\it induced subgraph of $G$ on $W$}, denoted by $G[W]$, is the graph on the vertex set $W$ with the edge
        set $\{e \in E(G) ~|~ e \subseteq W\}$. Note that $G[W]=G\setminus (V(G)\setminus W)$.

    \end{enumerate}

\noindent   
\textbf{Various classes of simple graphs:} 
\begin{enumerate}
    \item A {\it path} graph $P_n$ of length $n$ is a graph on the vertex set $\{x_1,\ldots,x_{n+1}\}$ with the edge set $\{\{x_i,x_{i+1}\}\mid 1\le i\le n\}$. A {\it cycle} $C_n$ of length $n$ is a graph on the vertex set $\{x_1,\ldots,x_n\}$ with the edge set $\{\{x_1,x_n\},\{x_i,x_{i+1}\}\mid 1\le i\le n-1\}$.

    \item A \textit{forest} is a graph without any induced cycle, and a connected forest is called a {\it tree}.

    \item The {\it star graph} $S_m$ of order $m$ is a tree with $m$ vertices among which one vertex has degree $m-1$, while the other $m-1$ vertices each has degree $1$.

    \item For $m>0$, a \textit{complete graph} $K_{m}$ is a graph on $m$ vertices such that there is an edge between any two distinct vertices.

    \item A graph $G$ is called \textit{chordal} if it contains no induced $C_n$ for $n\ge 4$. If $G$ is a chordal graph, then $G$ contains at least one vertex $y$ such that $N_G(y)$ is a complete graph (see \cite{Dirac}). We call such a vertex a {\it simplicial vertex} of $G$. Note that an induced subgraph of a chordal graph is again a chordal graph.
    \end{enumerate}

\subsection{The t-path ideal:} 
Let $G$ be a graph on the vertex set $\{x_1,\ldots,x_n\}$ and let $R$ denote the polynomial ring $\mathbb K[x_1,\ldots,x_n]$, where $\mathbb K$ is a field. A {\it $t$-path} from a vertex $x_p$ in $G$ to another vertex $x_q$ is a sequence of vertices $x_p=x_{i_1},x_{i_2},\ldots,x_{i_t}=x_q$ such that $\{x_{i_j},x_{i_{j+1}}\}\in E(G)$ for each $j\in [t-1]$.
 The $t$-path ideal of $G$, first introduced in \cite{ConcaDeNegri1998} and denoted by $I_{t}(G)$, is a square-free monomial ideal of $R$ given by
\[
I_t(G)=\left\l \left\{\prod_{j=1}^tx_{i_j}\mid x_{i_1},x_{i_2},\ldots,x_{i_t} \text{ forms a $t$-path in }G \right\}\right\r.
\]

\subsection{Path ideal as edge ideal of a hypergraph:}
A (simple) {\it hypergraph} $\H$ is a pair $(V(\H),E(\H))$, where $E(\H)\subseteq 2^{V(\H)}$ and for any two $\E,\E'\in E(\H)$, $\E\not\subset \E'$. The sets $V(\H)$ and $E(\H)$ are called the vertex set and edge set of $\H$, respectively. For a fixed positive integer $r$, if $|\E|=r$ for all $\E\in E(\H)$, then we say that $\H$ is an {\it $r$-uniform hypergraph}. Note that if $|\E|=2$ for each $\E\in E(\H)$, then $\H$ is just a graph. Thus, graphs are simply two uniform hypergraphs. If $\H$ is a hypergraph and $A\subseteq V(\H)$, then $\H\setminus A$ denotes the hypergraph on the vertex set $V(\H)\setminus A$ with the edge set $\{\E\in E(\H)\mid \E\cap A=\emptyset\}$. We say $\H'=\H\setminus A$ for some $A\subseteq V(\H)$ an \textit{induced subhypergraph} of $\H$ on the vertex set $V(\H)\setminus A$. For a vertex $x\in V(\H)$, we simply write $\H\setminus x$ to denote $\H\setminus\{x\}$. A \textit{vertex cover} of $\H$ is a collection of vertices $U\subseteq V(\H)$ such that for any edge $\E\in E(\H)$ one has $\E\cap U\neq \emptyset$. A \textit{minimal vertex cover} is a vertex cover which is minimal with respect to inclusion. 

Let $\H$ be a hypergraph on the vertex set $\{x_1,\ldots,x_n\}$ and let $R=\K[x_1,\ldots,x_n]$. Corresponding to each $\E\in E(\H)$, we assign the monomial $\x_{\E}=\prod_{x_j\in \E}x_j$ in $R$. Then the ideal $\l \x_{\E}\mid \E\in E(\H) \r$ is called the {\it edge ideal} of $\H$, and is denoted by $I(\H)$. Let $I\subseteq R$ be a square-free monomial ideal with the unique minimal monomial generating set $\G(I)$. Then $I$ can be viewed as an edge ideal of a hypergraph $\H_{I}$ by taking $V(\H_{I})=\{x_1,\ldots,x_n\}$ and $E(\H_{I})=\{\{x_{i_1},\ldots,x_{i_r}\}\mid x_{i_1}\cdots x_{i_r}\in \G(I)\}$, in other words, $I=I(\H_{I})$. It is well-known in the literature that the minimal prime ideals of $I$ (equivalently, the associated primes of $I$ as $I$ is radical) are exactly the ideals generated by the minimal vertex covers of $\H_{I}$. Therefore, the \textit{height} of $I$ (resp., the \textit{big height} of $I$), denoted by $\h(I)$ (resp., $\b(I)$), is the minimum (resp., maximum) cardinality of a minimal vertex cover of $\H_{I}$.\par

Let $G$ be a graph on the vertex set $\{x_1,\ldots,x_n\}$. Consider the ideal $I_{t}(G)$. Since $I_{t}(G)$ is a square-free monomial ideal, from the previous discussion, we can associate a hypergraph, say $\P_t(G)$, on the vertex set $\{x_1,\ldots,x_n\}$ such that $I_t(G)=I(\P_{t}(G))$. More precisely,
\begin{enumerate}
    \item[$\bullet$] $V(\P_t(G))=V(G)$,
    \item[$\bullet$] $E(\P_t(G))=\{\{x_{i_1},\ldots,x_{i_t}\}\mid x_{i_1},\ldots,x_{i_t} \text{ forms a } t\text{-path in } G\}$.
\end{enumerate}

\subsection{Some algebraic invariants:}
Let $I$ be a graded ideal in the polynomial ring $R=\K[x_1,\ldots , x_n]$. Then a graded minimal free resolution of $R/I$ is an exact sequence
\[
\mathcal F_{\cdot}: \,\, 0\rightarrow F_k\xrightarrow{\partial_{k}}\cdots\xrightarrow{\partial_{2}} F_1\xrightarrow{\partial_1} F_0\xrightarrow{\theta} R/I\rightarrow 0, 
\]
where $F_0=R$, $F_i=\oplus_{j\in\mathbb N}R(-j)^{\beta_{i,j}(R/I)}$ for $i\ge 1$, $\theta$ is the natural quotient map, and $R(-j)$ is the polynomial ring $R$ with its grading twisted by $j$. The numbers $\beta_{i,j}(R/I)$ are uniquely determined, and called the $i^{th}$ $\mathbb N$-graded {\it Betti numbers} of $R/I$ in degree $j$. The {\it Castelnuovo-Mumford regularity} of $R/I$, denoted by $\reg(R/I)$, is the number $\max\{j-i\mid \beta_{i,j}(R/I)\neq 0\}$. The invariant $\max\{i\mid \beta_{i,j}(R/I)\neq 0\}$ is called the {\it projective dimension} of $R/I$, and is denoted by $\pd(R/I)$. 
\par 

The following are some well-known results regarding regularity and projective dimension, which will be used throughout this paper.

\begin{lemma}\label{reg sum}\cite{HHBook}
        Let $I_1\subseteq R_1=\mathbb K[x_1,\ldots,x_n]$ and $I_2\subseteq R_2=\mathbb K[y_1,\ldots, y_m]$ be two graded ideals. Consider the ideal $I=I_1R+I_2R\subseteq R=\mathbb K[x_1,\ldots,x_n,y_1,\ldots,y_m]$. Then 
        \[ \reg(R/I)=\reg(R_1/I_1)+\reg(R_2/I_2)\,\,\text{and}\,\,\pd(R/I)=\pd(R_1/I_1)+\pd(R_2/I_2).
        \]
    \end{lemma}

\begin{lemma}\cite[Lemma 2.10, Lemma 5.1]{DHS}\label{regularity lemma}
    Let $I\subseteq R$ be a square-free monomial ideal and let $x_i$ be a variable appearing in some generator of $I$. Then
    \begin{enumerate}
        \item $   \reg(R/I)\le\max\{\reg(R/(I:x_i))+1,\reg(R/\langle I,x_i\rangle)\},$
        \item $\pd(R/I)\le\max\{\pd(R/(I:x_i)),\pd(R/\langle I,x_i\rangle)\}.$
    \end{enumerate}
    \end{lemma}

\begin{lemma}\label{regularity lemma1}\textup{(cf. \cite[Chapter 18]{IPBook})}
    Let $J$ and $K$ be two graded ideals of $R$. Then 
    \begin{enumerate}
        \item $\reg(R/(J+K))\le \max\{\reg(R/J),\reg(R/K),\reg(R/(J\cap K))-1\},$
        \item $  \pd(R/(J+K))\le \max\{\pd(R/J),\pd(R/K),\pd(R/(J\cap K))+1\}.$
    \end{enumerate}
\end{lemma}

\subsection{Bounds on regularity and projective dimension:}\label{induced_matching}
Let $\H$ be a hypergraph. A \textit{matching} in $\H$ is a collection of pairwise disjoint edges, that is, a subset $S\subseteq E(\H)$ such that for any two edges $\E,\E'\in S$ with $\E\neq \E'$, one has $\E\cap \E'=\emptyset$. An \textit{induced matching} is a matching in $\H$ for which the induced subhypergraph of $\H$ on the vertices involving the matching is a disjoint collection of edges. The following lower bound on the regularity in terms of the induced matching is well-known.

\begin{lemma}\cite[Theorem 4.2]{Ha2014}
    Let $\H$ be a hypergraph and $M$ be an induced matching of $\H$. Then \[\reg(R/I(\H))\geq \sum_{\E\in M}(|\E|-1).\]
\end{lemma}

\noindent

Let us define $\nu(\H)=\max\{|M|:M \text{ is an induced matching of }\H\}$, and call this to be the \textit{induced matching number} of the hypergraph $\H$. Then for a simple graph $G$, $\nu(G)$ gives a crude lower bound of $\reg(R/I(G))$. For our purpose, given a simple graph $G$, we call an induced matching of the hypergraph $\P_t(G)$ to be a \textit{$t$-path induced matching of $G$}. This is nothing but an induced subgraph of $G$ consisting of pairwise disjoint $t$-paths. The induced matching number of $\P_t(G)$ is denoted by $\nu_t(G)$ and is called the \textit{$t$-path induced matching number} of $G$. Therefore, in our case, we have the following lower bound for the regularity of $t$-path ideals.
\begin{lemma}\label{lower bound}
    Let $G$ be a finite simple graph. Then $\reg(R/I_t(G))\ge (t-1)\nu_t(G)$.
\end{lemma}

Using the Alexander dual of square-free monomial ideals (see Section \ref{section 4}) and Terai's formula \cite[Theorem 2.1]{Terai}, one  can get an analogous bound for the projective dimension of $R/I_t(G)$ in terms of $\b(I_t(G))$ as follows:

\begin{lemma}\label{lempdbght}
    Let $G$ be a finite simple graph. Then $\pd(R/I_t(G))\ge \b(I_t(G))$.
\end{lemma}

\section{\texorpdfstring{$3$-path}{} ideals of chordal graphs}\label{section 3}

In this section, our goal is to find the regularity and projective dimension of $3$-path ideals of chordal graphs. In particular, we use $\nu_3(G)$ to express the regularity and $\b(I(\P_3(G)))$ to express the projective dimension of $R/I_{3}(G)$ for any chordal graph $G$. We use a new splitting of ideals to get our desired results. First, let us start by deriving some results needed for the proof of our main theorems. \par 

A chordal graph $G$ always contains a simplicial vertex, i.e. a vertex $x\in V(G)$ such that the induced subgraph on the set of vertices $N_G[x]$ is a complete graph. As a consequence, for each $y\in N_G(x)$, one has $N_G[x]\subseteq N_G[y]$. In general, for any graph $G$, if there are two vertices $x, y\in V(G)$ such that $N_G[x]\subseteq N_G[y]$, then $\nu_3(G)$ and $\nu_3(G- e)$, where $e=\{x,y\}$, are related in the following way.

\begin{proposition}\label{3 path matching}
    Let $G$ be a graph and $x,y\in V(G)$ such that $N_G[x]\subseteq N_G[y]$. Then $\nu_3(G- e)\le\nu_3(G) $, where $e=\{x,y\}$.
\end{proposition}
\begin{proof}
    If $x\in V(G)$ is a leaf vertex, then the assertion is true. So, we assume that $x$ is not a leaf vertex. Let $A$ be a maximal $3$-path induced matching of $G- e$ such that $\nu_3(G- e)=|A|$. Our aim is to show that $A$ is also a $3$-path induced matching of $G$. First assume that for each $3$-path $\mathcal E\in A$ of $G- e$, $x\notin \mathcal E$. Then it is easy to see that $A$ is also a $3$-path induced matching of $G$. Consequently, $\nu_3(G- e)\le \nu_3(G)$. Now suppose that $\mathcal E\in A$ is a $3$-path of $G- e$ such that $x\in\mathcal E$. Observe that if there is a 3-path $\mathcal{E}'\in A, \mathcal{E}'\neq \mathcal{E}$ such that $y\in \mathcal{E}'$, then $A$ does not give a $3$-path induced matching of $G$. We shall show that this is not the case. Since $x\in \mathcal{E}$, there exists some $z\in (N_G(x)\setminus \{y\})\subseteq (N_G(y)\setminus \{x\})$ such that $z\in\mathcal E$. Hence, for any $\mathcal E'\in A$ and $\mathcal E\neq \mathcal E'$, we have $y\notin \mathcal E'$, as $y\in N_{G-e}(z)$. Thus, there exists at most one $3$-path $\mathcal E\in A$ of $G-e$ such that $\mathcal E\cap\{x,y\}\neq \emptyset$. This implies $A$ is also a $3$-path induced matching of $G$. Hence, $\nu_3(G-e)\le \nu_3(G)$.
\end{proof}

\begin{corollary}\label{matching result}
    Let $G$ be a graph and $x\in V(G)$ be a simplicial vertex of $G$. Then for any $y\in N_G(x)$, we have $\nu_3(G- e)\le \nu_3(G)$, where $e=\{x,y\}$.
\end{corollary}

\begin{lemma}\label{1 less induced matching}
    Let $G$ be a graph and $x,y\in V(G)$ such that $N_G[x]\subseteq N_G[y]$ and $N_G(y)=\{x,w_1,\ldots,w_r\}$. Then $\nu_3(G\setminus N_G[y,w_i])\le \nu_3(G)-1$ for each $i\in [r]$.
\end{lemma}
\begin{proof}
    Let $A$ be a $3$-path induced matching of $G\setminus N_G[y,w_i]$ such that $|A|=\nu_3(G\setminus N_G[y,w_i])$. If $\E\in A$, then it is clear that for each $z\in \E$, $\{z,y\}\notin E(G)$ and $\{z,w_i\}\notin E(G)$. Moreover, since $\{x,y\}\in E(G)$, we also have $\{x,z\}\notin E(G)$. Therefore, $A\cup\{x,y,w_i\}$ is a $3$-path induced matching of $G$. Hence, $\nu_3(G\setminus N_G[y,z])\le \nu_3(G)-1$.
\end{proof}
The following are some easy observations that will be used repeatedly.

\begin{proposition}\label{connected comp lemma}
     Let $G=G_1\sqcup G_2$ be a disconnected graph. Then $\nu_3(G)=\nu_3(G_1)+\nu_3(G_2)$.
\end{proposition}

\begin{proposition}\label{G-e chordal}
    Let $G$ be a chordal graph and $x$ be a simplicial vertex of $G$. Then for any $y\in N_{G}(x)$, $G-e$ is chordal, where $e=\{x,y\}$.
\end{proposition}
It is important to note that if $G$ is a chordal graph and $e\in E(G)$ is an arbitrary edge, then $G-e$ need not be chordal.

\begin{proposition}\label{colon comma exchange}
    Let $I\subseteq R$ be a monomial ideal, and $x_1,\ldots,x_r\in R$ are some indeterminates. Then \[ (I:x_r)+\l x_1,\ldots,x_{r-1} \r=(( I+\l x_1,\ldots,x_{r-1} \r):x_r).\]
\end{proposition}

    To prove the main result, our aim is to write the $3$-path ideal $I_3(G)=J+K$, for some suitable choice of $J$ and $K$ so that $\reg(R/J)$, $\reg(R/K)$, and $\reg(R/(J\cap K))$ are relatively easier to determine. For this, we make use of the fact that $G$ is chordal, and hence, $G$ contains a simplicial vertex, say $x\in V(G)$. We split the ideal $I_3(G)$ by means of any edge containing the simplicial vertex $x$. Then, as we shall see in our main theorem, using this decomposition of the ideal $I_3(G)$, we get the desired bounds for $\reg(R/J)$, $\reg(R/K)$, and $\reg(R/(J\cap K))$. This new splitting of the $3$-path ideal has not been used in the literature much as there are some generators of $I_3(G)$ that appear in both $J$ and $K$. Researchers usually tend to avoid such a situation due to the complexity of the arguments that it might add. But, for our case, this allows us to identify $J, K, \text{ and }J\cap K$ with some smaller combinatorial objects and calculate the regularity and projective dimension by using some inductive arguments. In the following lemma, we introduce the splitting of $I_3(G)$ as mentioned above, which will be crucial in the proofs of \Cref{result1} \& \ref{pdmain}.

 \begin{lemma}\label{auxiliary lemma}
    Let $x$ and $y$ be two vertices of a graph $G$ such that $N_G[x]\subseteq N_G[y]$. Consider the ideals $J=\l xyw\mid w\in N_G(y)\setminus \{x\} \r \text{ and } K=I_3(G- e),$ where $e=\{x,y\}$. Then  
    \begin{enumerate}
        \item $I_3(G)=J+K$,
        \item $J\cap K=xyL$, where
        \begin{itemize}
            \item[(a)] $(L:w)=\l\{ v\mid v\in N_G[w,y]\setminus\{x,y,w\}\}\r+I_3(G\setminus N_G[w,y])$ for all $w\in N_G(y)\setminus N_G[x],$ 
            \item[(b)] $(L:w)=R$ for all $w\in N_G(x)\setminus \{y\}$. 
        \end{itemize}
    \end{enumerate}
\end{lemma}
\begin{proof}
$(1)$ is immediate. We proceed to prove $(2)$.
Let $N_G(y)=\{x,w_1,\ldots,w_r,w_{r+1},\ldots,w_{r+s}\}$, where $\{ x,w_i\}\notin E(G)$ for each $i\in [r]$ and $\{x,w_i\}\in E(G)$ for each $r+1\le i\le r+s$. In other words, $N_G(y)\setminus N_G[x]=\{w_1,\ldots,w_r\}$ and $N_G(x)\setminus \{y\}=\{w_{r+1},\ldots,w_{r+s}\}$. Then we have $J=\l xyw_1,\ldots,xyw_r,xyw_{r+1},\ldots,xyw_{r+s}\r$. Also, $xyw_i\notin \G(K)$ for all $i\in[r]$ and $xyw_i\in \G(K)$ for all $r+1\le i\le r+s$. It is easy to see that $xy$ divides each of the monomials $\lcm(m,m')$, where $m\in J$ and $m\in K$. Thus, $J\cap K=xyL$, where $L$ is generated by the monomials $\frac{\lcm(m,m')}{xy}$ for $m\in J$ and $m'\in K$. Since $xyw_i\in \G(J\cap K)$ for each $r+1\le i\le r+s $, we have $(L:w_i)=R$ for such an $i$. In other words, $(L:w)=R$ for all $w\in N_G(x)\setminus\{y\}$, which proves $(b)$.

To prove $(a)$, we need to show that \[(L:w_i)=\l\{ v\mid v\in N_G[w_i,y]\setminus\{x,y,w_i\}\}\r+I_3(G\setminus N_G[w_i,y])\] for each $i\in[r]$. Without any loss of generality, let us assume that $i=1$. Consider the ideal \[M=\l\{ v\mid v\in N_G[w_{1},y]\setminus\{x,y,w_{1}\}\}\r+I_3(G\setminus N_G[w_{1},y]).\] Then for each $v\in N_G(w_{1})\setminus \{y\}$, $xyw_{1}v=\lcm (xyw_{1},yw_{1}v)$, where $xyw_{1}\in \G(J)$ and $yw_{1}v\in \G(K)$. Also, if $v\in N_G(y)\setminus \{x,w_{1}\}$, then $v\in \{w_2,\ldots,w_r,w_{r+1},\ldots,w_{r+s}\}$. Hence, $xyw_{1}v=\lcm(xyv,vyw_{1})$, where $xyv\in \G(J)$ and $vyw_{1}\in \G(K)$. Note that $I_3(G\setminus N_G[w_{1},y])\subseteq K$. Thus, if $abc\in\G(I_3(G\setminus N_G[w_{1},y]))$, then $xyw_{1}abc=\lcm (xyw_{1},abc)$, where $xyw_{1}\in \G(J)$ and $abc\in \G(K)$. Consequently, $M\subseteq (L:w_{1})$.

To prove the other containment, we first determine the generators of the ideal $J\cap K$. Let $N_G(w_{1})=\{y,p_1,\ldots,p_k\}$. It may happen that $p_i=w_j$ for some $j\ge 2$. Now the ideal $J\cap K$ is generated by the monomials $\lcm(m,m')$, where $m'\in I_3(G- e)$ and $m=xyw_i$ for some $i\ge 1$. Set
    \begin{align*}
        A_1&=\{\lcm(m,m')\mid w_i|m\text{ or }w_i|m' \text{ for some }i\ge 2,\text{ where }m\in J,m'\in K \},\\
        A_2&=\{\lcm(m,m')\mid p_t|m \text{ or } p_t|m' \text{ for some }t\ge 1,\text{ where }m\in J,m'\in K \},\\
        A_3&=\{\lcm(xyw_{1},abc)\mid abc\in I_3(G\setminus N_G[w_{1},y])\}.
    \end{align*}
      Then $J\cap K$ is generated by the monomials in $\cup_{i=1}^3A_i$. It is easy to see that if $\widetilde m\in A_1$, then $xyw_i$ divides $\widetilde m$ for some $i\ge 2$. Similarly, if $\widetilde m\in A_2$ then $xyp_t$ divides $\widetilde m$ for some $t\ge  1$. Also, $A_3=\{xyw_{1}abc\mid abc\in I_3(G\setminus N_G[y,w_{1}])\}$. Note that $xyw_{1}w_i\in A_1$ for each $i\ge 2$, and $xyw_{1}p_t\in A_2$ for each $t\ge 1$. Using this description of $J\cap K$ we conclude that $(L:w_{1})\subseteq M$. This completes the proof.
\end{proof}

Now, we prove one of the main theorems of this section, which provides a combinatorial formula for the regularity of $3$-path ideals of chordal graphs. The key idea of the proof is the splitting of the ideal $I_3(G)$ as mentioned in \Cref{auxiliary lemma}.

\begin{theorem}\label{result1}
    Let $G$ be a chordal graph. Then
    \[
    \reg(R/I_3(G))= 2\nu_3(G),
    \]
where $\nu_3(G)$ denotes the $3$-path induced matching number of $G$.
\end{theorem}
\begin{proof}
    In view of \cref{lower bound}, we only need to show that $\reg(R/I_3(G))\le 2\nu_3(G)$. We proceed by induction on $|E(G)|$. If $|E(G)|=1$, then $\reg(R/I_3(G))=0=2\nu_3(G)$. Now suppose $|E(G)|=2$. Then either $I_3(G)=\l 0\r$ or $I_3(G)=\l abc\r$, for some $a,b,c\in V(G)$. In both cases, it is easy to see that $\reg(R/I_3(G))=2\nu_3(G)$. Therefore, we may assume that $|E(G)|\ge 3$.

        First, let us consider the case when $G$ is a complete graph. Then $I_3(G)$ is an edge ideal of a $3$-uniform complete hypergraph. Then by \cite[Theorem 3.1]{Emtander2009}, $I_3(G)$ has a linear resolution. Therefore,  $\reg(R/I_3(G))=2=2\nu_3(G)$. Hence, due to \Cref{connected comp lemma} and \Cref{reg sum}, we may assume that $G$ is not a disjoint union of complete graphs. In that case, there exists some $x,y\in V(G)$ such that $x$ is a simplicial vertex of $G$, $y\in N_G(x)$, and $N_G[x]\subsetneq N_G[y]$. Without loss of generality, let $N_G(y)=\{x,w_1,\ldots,w_r,w_{r+1},\ldots,w_{r+s}\}$, where $\{x,w_i\}\notin E(G)$ for $i\in [r]$ and $\{x,w_i\}\in E(G)$ for $r+1\le i\le r+s$. Note that $r$ is a positive integer as $N_G[x]\subsetneq N_G[y]$. Now, let us consider the ideals
    \[J=\l xyw\mid w\in N_G(y)\setminus \{x\}\r \text{ and }K=I_3(G- e),
    \]
    where $e=\{x,y\}\in E(G)$. Then $J+K=I_3(G)$, $J=\l xyw_1,\ldots,xyw_{r+s}\r$, $xyw_i\in K$ for $r+1\le i\le r+s$ and $xyw_i\notin K$ for each $i\in [r]$.  Moreover, $J\cap K=xyL$, where 
    \[
    L=\left\l\left\{ \frac{\lcm(m,m')}{xy}\mid m\in J\text{ and }m'\in K\right\}\right\r.
    \]
    By the construction of $J$ and $K$, we have $w_i\notin \G(L)$ for $i\in [r]$ and $w_i\in \G(L)$ for $r+1\le i\le r+s$. Based on these observations, we prove the following.
    
    \noindent
    {\bf Claim 1}: For each $i\in [r]$, $\reg(R/(L:w_i))\le 2\nu_3(G)-2$.
    
    \noindent
    {\bf Proof of Claim 1}: Fix some $i\in [r]$. Let $N_G(w_i)\setminus N_G[y]=\{u_{i1},\ldots,u_{ik_i}\}$, where $k_i\ge 0$. Then by \Cref{auxiliary lemma},
    \[
    (L:w_i)=\l w_1,\ldots\widehat{w_i},\ldots,w_{r+s}\r+\l u_{i1},\ldots,u_{ik_i}\r+I_3(G\setminus N_G[y,w_i]).
    \] 

\begin{figure}
    \centering
    \begin{tikzpicture}
        [scale=.55]
        \draw [fill] (0,0) circle [radius=0.1];
        \draw [fill] (0,2.5) circle [radius=0.1];
        \draw [fill] (2.5,0) circle [radius=0.1];
        \draw [fill] (2.5,2.5) circle [radius=0.1];
        \draw [fill] (5,5) circle [radius=0.1];
        \draw [fill] (5,1) circle [radius=0.1];
        \draw [fill] (7.5,4) circle [radius=0.1];
        \draw [fill] (5,3) circle [radius=0.1];
        \draw [fill] (7.5,2) circle [radius=0.1];
        \node at (-0.5,0) {};
        \node at (0,-0.5) {$w_{r+s}$};
        \node at (0,3) {$x$};
        \node at (2.5,-0.5) {$w_{r+1}$};
        \node at (2.5,3) {$y$};
        \node at (5,5.5) {$w_1$};
        \node at (5,0.5) {$w_r$};
        \node at (5,2.5) {$w_{i}$};
        \node at (7.5,4.5) {$u_{i,1}$};
        \node at (7.5,1.5) {$u_{i,k_i}$};
        \node at (1.25,0) {$\cdots$};
        \node at (5,4) {$\vdots$};
        \node at (5,2) {$\vdots$};
        \node at (7.5,3.2) {$\vdots$};
        \draw [thick,red] (0,2.5)--(2.5,2.5);
        \draw (0,2.5)--(0,0);
        \draw (2.5,2.5)--(0,0);
        \draw (2.5,2.5)--(5,5);
        \draw (2.5,2.5)--(2.5,0)--(0,2.5);
        \draw (2.5,2.5)--(5,1);
        \draw (2.5,2.5)--(5,3);
        \draw (5,3)--(7.5,2);
        \draw (5,3)--(7.5,4);
    \end{tikzpicture}
    \caption{Splitting of a chordal graph $G$ by means of the edge $\{x,y\}$.}
    \label{fig:enter-label0}
\end{figure}
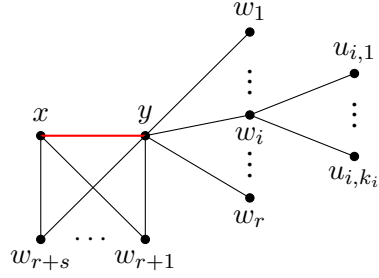

    \noindent
    Therefore, by the induction hypothesis, $\reg(R/(L:w_i))=\reg(R/I_3(G\setminus N_G[y,w_i]))\le 2\nu_3(G\setminus N_G[y,w_i])$. Now by \Cref{1 less induced matching}, we have $\nu_3(G\setminus N_G[y,w_i]))\le \nu_3(G)-1$. Consequently, $\reg(R/(L:w_i))\le 2\nu_3(G)-2$.
    
    \noindent
    {\bf Claim 2}: $\reg(R/L)\le 2\nu_3(G)-1$.
    
    \noindent
    {\bf Proof of Claim 2}: Note that $L+\l w_1,\ldots,w_r \r=\l w_1,\ldots,w_{r+s} \r$, and this implies $\reg(R/(L+\l w_1,\ldots,w_r \r))=0$. By our assumption, $\nu_3(G)\ge 1$ as $\{x,y,w_1\}$ is a $3$-path induced matching of $G$. Consequently, $\reg(R/(L+\l w_1,\ldots,w_r \r))\le 2\nu_3(G)-1$. Now consider the ideal $(L+\l w_1,\ldots,w_{r-1}\r):w_r$. By \Cref{colon comma exchange}, $(L+\l w_1,\ldots,w_{r-1}\r):w_r=(L:w_r)+\l w_1,\ldots,w_{r-1} \r$. Therefore, using \Cref{auxiliary lemma}, we obtain $(L+\l w_1,\ldots,w_{r-1}\r):w_r=(L:w_r)$, and hence, by Claim 1, $\reg(R/(L+\l w_1,\ldots,w_{r-1}\r):w_r)\le  2\nu_3(G)-2$. Thus, using \Cref{regularity lemma}, we obtain $\reg(R/(L+\l w_1,\ldots,w_{r-1} \r))\le 2\nu_3(G)-1$. Next we consider the ideal $(L+\l w_1,\ldots,w_{r-2} \r):w_{r-1}$ and proceed as before. By \Cref{colon comma exchange},  $(L+\l w_1,\ldots,w_{r-2} \r):w_{r-1}=(L:w_{r-1})+\l w_1,\ldots,w_{r-2} \r$. Then using \Cref{auxiliary lemma} again, we get $(L+\l w_1,\ldots,w_{r-2} \r):w_{r-1}=(L:w_{r-1})$, and hence, by Claim 1, $\reg(R/(L+\l w_1,\ldots,w_{r-2} \r):w_{r-1})\le 2\nu_3(G)-2$. Therefore, again using \Cref{regularity lemma}, we obtain $\reg(R/(L+\l w_1,\ldots,w_{r-2} \r))\le 2\nu_3(G)-1$. Continuing in this way after a finite number of steps, we get $\reg(R/L)\le 2\nu_3(G)-1$.

Now, since $J\cap K=xyL$ and both $x$ and $y$ do not divide any minimal generator of $L$, we conclude that $\reg(R/J\cap K)\le 2\nu_3(G)+1$. For the ideal $J$, we have $\reg(R/J)=2\le 2\nu_3(G)$ as $\nu_3(G)\ge 1$. Also, using \Cref{G-e chordal} and by the induction hypothesis, $\reg(R/K)=\reg(R/I_3(G-e))\le 2\nu_3(G- e)$, where $e=\{x,y\}$. Moreover, by \Cref{matching result}, $\nu_3(G-e)\le \nu_3(G)$, and thus, $\reg(R/K)\le 2\nu_3(G)$. Therefore, using \Cref{regularity lemma1}, we conclude that $\reg(R/I_3(G))\le 2\nu_3(G)$. This completes the proof of the theorem.
\end{proof}

We now move on to find the projective dimension of $3$-path ideals of chordal graphs. For this, our approach is again finding a suitable decomposition of the $3$-path ideals that is compatible with the inductive arguments. We shall see that the same decomposition as in the case of regularity works for the case of projective dimension as well. Specifically, we prove that the projective dimension of the $3$-path ideal of a chordal graph is equal to the big height of that ideal. Before going into the proof, we present a few preparatory lemmas related to the big height of $3$-path ideals. 

\begin{lemma}\label{connected comp pd}
    Let $G=G_1\sqcup G_2$ be a disconnected graph. Then $\b(I_3(G))=\b(I_3(G_1))+\b(I_3(G_2))$.
\end{lemma}
\begin{proof}
    Follows from the fact that for a square-free monomial ideal $I=I(\H)$, $\b(I)$ is the maximum cardinality of a minimal vertex cover of the hypergraph $\H$.
\end{proof}

\begin{lemma}\label{induced bight}
    Let $x,y\in V(G)$ such that $N_G[x]\subseteq N_G[y]$. Then $\b(I_3(G-e))\le\b(I_3(G))$, where $e=\{x,y\}$.
\end{lemma}
\begin{proof}
    Let $H=G-e$ and let $\p$ be a minimal prime ideal of $I_3(H)$ such that $\h(\p)=\b(I_3(H))$. If $x\in \p$ or $y\in\p$, then it is easy to see that $\p$ is a minimal prime ideal of $I_3(G)$. Hence, $\b(I_3(H))\le\b(I_3(G))$. Therefore, we may assume that $x,y\notin \p$. Now let $N_G(y)=\{x,w_1,\ldots,w_s\}$. If $w_i\in\p$ for all $i\in [s]$, then also $\p$ is a minimal prime ideal of $I_3(G)$, and hence, $\b(I_3(H))\le\b(I_3(G))$. Without loss of generality, let $w_1\notin\p$. Then $\{x,w_1\}\notin E(H)$ and $w_i\in\p$ for all $2\le i\le s$. We proceed to show that $\p+\l x\r$ is a minimal prime ideal of $I_3(G)$. Since $H=G- e$, we see that $I_3(G)\subseteq \p+\l x\r$. Let $\q$ be the minimal prime ideal of $I_3(G)$ such that $\q\subseteq \p+\l x\r$. Since $xyw_1\in I_3(G)$ and $y,w_1\not\in \p$, we should have $x\in\q$. If $\q\subsetneq\p+\l x\r$, then there exists some $z\in\G(\p)$ such that $z\notin \G(\q)$. Note that $z\neq x$ as $x\in \q$, and $z\neq w_i$ for each $2\le i\le s$ as $w_1yw_i\in I_3(G)$ and $w_1,y\notin \q$. Now, let $\p'=\l u\in \G(\p)\mid u\neq z\r$ and $abc\in \G(I_3(H))$. Observe that $\p'\subsetneq\p$. We want to show that $abc\in\p'$. If $x\nmid abc$, then $abc\in\p'$ since $abc\in\q$. Now, suppose $x\mid abc$. Then $w_i\mid abc$ for some $i\in [s]$ since $N_G[x]\subseteq N_G[y]$. If $w_j\mid abc$ for some $2\le j\le s$, then clearly, $abc\in \p'$. Now if $w_1\mid abc$, then also $w_j\mid abc$ for some $2\le j\le s$ since $\{x,w_1\},\{x,y\}\notin E(H)$. Therefore, $abc\in\p'$ and thus, $I_3(H)\subseteq \p'$, a contradiction. Therefore, $\q=\p+\l x\r$ and hence, $\b(I_3(H))\le \b(I_3(G))$.   
\end{proof}

\begin{lemma}\label{bight prop}
    Let $G$ be a graph and $y\in V(G)$ with $N_G(y)=\{w_1,\ldots,w_r\}$. For $i\in[r]$, if $|N_G(w_i)\setminus N_G[y]|=k_i$, then $\b(I_3(G))\ge \b(I_3(G\setminus N_G[y,w_i]))+r+k_i-1$.
\end{lemma}

\begin{proof}
    Without loss of generality, take $i=1$. Let $\p$ be a minimal prime ideal of $I_3(G\setminus N_G[y,w_1])$ such that $\h(\p)=\b(I_3(G\setminus N_G[y,w_1]))$. Let $N_G(w_1)\setminus N_G[y]=\{u_{11},\ldots,u_{1k_1}\}$. Then it is easy to see that $\q=\p+\l w_2,\ldots,w_{r}\r+\l u_{11},\ldots,u_{1k_1}\r$ is a prime ideal containing $I_3(G)$. Now, suppose $\q$ is not minimal. Let $\q'\subsetneq\q$ be a minimal prime ideal of $I_3(G)$. Since $yw_1w_i\in I_3(G)$, we have $w_i\in\q'$ for each $2\le i\le r$. Similarly, since $yw_1u_{1j}\in I_3(G)$, we have $u_{1j}\in\q'$ for each $j\in [k_1]$. Hence, there exists some $z\in\G(\p)$ such that $z\notin\q'$. Let $\p'=\l u\in\G(\p)\mid u\neq z \r$. Note that if $abc\in\G(I_3(G\setminus N_G[y,w_1]))$, then $y\nmid abc,w_i\nmid abc,$ and $u_{1t}\nmid abc$ for each $i\in[r]$ and $t\in[k_1]$. Also, $abc\in I_3(G)$ implies $abc\in\q'$, and hence, $abc\in \p'$. Thus, $\p'\subsetneq \p$ is a prime ideal containing $I_3(G\setminus N_G[y,w_1])$, a contradiction. Therefore, $\q$ is a minimal prime ideal of $I_3(G)$. Note that $\h(\q)=\h(\p)+r+k_1-1$, and thus, the result follows.
\end{proof}

Now, we are ready to prove the second main result of this section.

\begin{theorem}\label{pdmain}
    Let $G$ be a chordal graph. Then
    \[
    \pd(R/I_3(G))=\b(I_3(G)),
    \]
    where $\b(I_3(G))$ denotes the maximum height of a minimal prime ideal of $I_3(G)$.
\end{theorem}

\begin{proof}
    
    In view of \Cref{lempdbght}, it is enough to prove that $\pd(R/I_3(G))\le\b(I_3(G))$, and we proceed to prove this inequality by induction on $|E(G)|$. If $|E(G)|=1$, then $\pd(R/I_3(G))=0=\b(I_3(G))$. Now suppose $|E(G)|=2$. Then either $I_3(G)=\l 0\r$ or $I_3(G)=\l abc\r$, for some $a,b,c\in V(G)$. In both cases, it is easy to see that $\pd(R/I_3(G))=\b(I_3(G))$. Therefore, we may assume that $|E(G)|\ge 3$.

    First, let us consider the case when $G$ is a complete graph. Then $I_3(G)^{\vee}$ is an edge ideal of a $(n-2)$-uniform complete hypergraph, and hence, by \cite[Theorem 3.1]{Emtander2009} along with the Terai's formula, we get $\pd(R/I_3(G))=n-2=\b(I_3(G))$. 
    Therefore, using \Cref{reg sum} and \Cref{connected comp pd} we may assume that $G$ is not a disjoint union of complete graphs. In that case, there exists some $x,y\in V(G)$ such that $x$ is a simplicial vertex of $G$, $y\in N_G(x)$, and $N_G[x]\subsetneq N_G[y]$. Without loss of generality, let $N_G(y)=\{x,w_1,\ldots,w_r,w_{r+1},\ldots,w_{r+s}\}$, where $\{x,w_i\}\notin E(G)$ for $i\in[r]$, and $\{x,w_i\}\in E(G)$ for $r+1\le i\le r+s$. Note here that $r$ is a positive integer since $N_G[x]\subsetneq N_G[y]$.
    Let us consider the ideals
    \[J=\l xyw\mid w\in N_G(y)\setminus \{x\}\r \text{ and }K=I_3(G- e),
    \]
    where $e=\{x,y\}\in E(G)$. Then $J+K=I_3(G)$, $J=\l xyw_1,\ldots,xyw_{r+s}\r$, $xyw_i\in K$ for $r+1\le i\le r+s$, and $xyw_i\notin K$ for each $i\in [r]$.
    
    If $V(G)=N_G[y]$, then $J\cap K=\l xyw_{r+1},\ldots ,xyw_{r+s}\r+\l xyw_iw_j\mid 1\le i<j\le r\r$ and hence, $\pd(R/J\cap K)=r+s-1$ . Note that $\pd(R/I_3(G))\ge r+s$ since $\l w_1,\ldots,w_{r+s}\r$ is a minimal prime ideal of $I_3(G)$. Thus, $\pd(R/J\cap K)\le\b(I_3(G))-1$. Also, $\pd(R/J)=r+s\le \b(I_3(G))$. Now using \Cref{G-e chordal} and by the induction hypothesis along with \Cref{induced bight}, $\pd(R/K)\le \b(I_3(G-e))\le\b(I_3(G))$. Therefore, using \Cref{regularity lemma1} we conclude that $\pd(R/I_3(G))\le \b(I_3(G))$.

    Now let $V(G)\neq N_G[y]$. Then either $I_3(G)=I_3(G[N_G[y]])$ or $I_3(G[N_G[y]])\subsetneq I_3(G)$. In the first case, proceeding as in the previous paragraph, we have $\pd(R/I_3(G))\le \b(I_3(G))$. Therefore, we may assume that $I_3(G[N_G[y]])\subsetneq I_3(G)$.

    \noindent
    {\bf Case I}: $G[N_G[y]]$ forms a connected component of $G$. In this case, using \Cref{reg sum} and \Cref{connected comp pd} together, we obtain $\pd(R/I_3(G))=\pd(R/I_3(G[N_G[y]]))+\pd(R/I_3(G[W]))$ and $\b(I_3(G))=\b(I_3(G[N_G[y]]))+\b(I_3(G[W]))$,
    where $W=V(G)\setminus N_G[y]$. Since $I_3(G[N_G[y]])\subsetneq I_3(G)$, we see that $E(G[W])\neq\emptyset$.  Therefore, using the induction hypothesis, we get $\pd(R/I_3(G[N_G[y]]))\le \b(I_3(G[N_G[y]]))$ and $\pd(R/I_3(G[W]))\le \b(I_3(G[W]))$. Consequently, $\pd(R/I_3(G))\le \b(I_3(G))$.

    \noindent
    {\bf Case II}: $G[N_G[y]]$ do not form a connected component of $G$. Recall that, $I_3(G)=J+K$, where $J=\l xyw_1,\ldots,xyw_{r+s}\r$ and $K=I_3(G-e)$, where $e=\{x,y\}\in E(G)$. Then $J\cap K=xyL$, where $w_i\notin \G(L)$ for $i\in [r]$, and $w_i\in\G(L)$ for $r+1\le i\le r+s$. Based on these observations, we prove the following.

    \noindent
    {\bf Claim 1}: $\b(I_3(G))\ge r+s+1$.

    \noindent
    {\bf Proof of claim 1}: Since $G[N_G[y]]$ do not form a connected component of $G$, there exists some $i\in [r+s]$ such that $N_G(w_i)\setminus N_G[y]\neq\emptyset$. For such an $i$, let $N_G(w_i)\setminus N_G[y]=\{u_{i1},\ldots,u_{ik_i}\}$, where $k_i\ge 1$. Then by \Cref{bight prop}, $\b(I_3(G))\ge \b(I_3(G\setminus N_G[y,w_i]))+r+s+k_i\ge r+s+1$.
    \noindent
    {\bf Claim 2}: For each $i\in [r]$, $\pd(R/(L:w_i))\le \b(I_3(G))-1$.

    \noindent
    {\bf Proof of Claim 2}: For $i\in [r]$, let $N_G(w_i)\setminus N_G[y]=\{u_{i1},\ldots,u_{ik_i}\}$, where $k_i\ge 0$. For such an $i$, by \Cref{auxiliary lemma}, 
    \[
    (L:w_i)=\l w_1,\ldots\widehat{w_i},\ldots,w_{r+s}\r+\l u_{i1},\ldots,u_{ik_i}\r+I_3(G\setminus N_G[y,w_i]).
    \]
    Therefore, by the induction hypothesis we have, $\pd(R/(L:w_i))\le r+s-1+k_i+\b(I_3(G\setminus N_G[y,w_i]))$. Note that $N_G(y)=\{x,w_1,\ldots,w_{r+s}\}$. Hence, by \Cref{bight prop}, $\b(I_3(G))\ge r+s+k_i+\b(I_3(G\setminus N_G[y,w_i]))$. Consequently, $\pd(R/(L:w_i))\le \b(I_3(G))-1$. 

    \noindent
    {\bf Claim 3}: $\pd(R/L)\le \b(I_3(G))-1$.

    \noindent
    {\bf Proof of Claim 3}: Since each monomial of $L$ is divisible by $w_i$ for some $i\in [r+s]$, we have $L+\l w_1,\ldots,w_r\r=\l w_1,\ldots,w_{r+s}\r$. Therefore, $\pd(R/L+\l w_1,\ldots,w_r\r)=r+s$, and hence, by Claim 1, we obtain $\pd(R/L+\l w_1,\ldots,w_r\r)\le \b(I_3(G))-1$. Now by \Cref{colon comma exchange}, $(L+\l w_1,\ldots, w_{r-1}\r):w_r=(L:w_r)+\l w_1,\ldots,w_{r-1}\r$. Therefore, using \Cref{auxiliary lemma}, we get $(L+\l w_1,\ldots, w_{r-1}\r):w_r=(L:w_r)$. Note that $\pd(R/(L:w_r))\le\b(I_3(G))-1$ because of Claim 2, and hence, $\pd(R/(L+\l w_1,\ldots, w_{r-1}\r):w_r)\le\b(I_3(G))-1$. Consequently, by \Cref{regularity lemma}, we have $\pd(R/L+\l w_1,\ldots,w_{r-1}\r)\le\b(I_3(G))-1$. Now $(L+\l w_1,\ldots,w_{r-2}\r):w_{r-1}=(L:w_{r-1})+\l w_1,\ldots,w_{r-2}\r$ (again by \Cref{colon comma exchange}). Therefore, using \Cref{auxiliary lemma}, Claim 2, and \Cref{regularity lemma} as above, we obtain $\pd(R/L+\l w_1,\ldots,w_{r-2}\r)\le\b(I_3(G))-1$. Continuing the above process, we finally obtain $\pd(R/L)\le \b(I_3(G))-1$ as desired. 

    Note that $J\cap K=xyL$, and hence, $\pd(R/J\cap K)\le \b(I_3(G))-1$. Now $\pd(R/J)=r+s<\b(I_3(G))$, and by using \Cref{induced bight}, \Cref{G-e chordal} and by the induction hypothesis, we have $\pd(R/K)\le\b(I_3(G-e))\le\b(I_3(G))$. Therefore, $\pd(R/I_3(G))\le \b(I_3(G))$ (by \Cref{regularity lemma1}). This completes the proof of the theorem.
\end{proof}

In \cite{CamposGundersonMoreyPaulsenPolstra2014}, the authors first studied Cohen-Macaulayness of path ideals (of any length), and characterized this property in the case of trees. As an immediate consequence of \Cref{pdmain}, we get the following equivalent criteria for Cohen-Macaulayness of $3$-path ideals of chordal graphs.
\begin{corollary}\label{CM}
    Let $G$ be a chordal graph. Then $I_3(G)$ is Cohen-Macaulay if and only if it is unmixed.
\end{corollary}
\begin{proof}
    Assume that $I_3(G)$ is unmixed. Then by \Cref{pdmain} and by the Auslander-Buchsbaum formula, 
    \[\d(R/I_3(G))=\dim(R)-\b(I_3(G))=\dim(R)-\h(I_3(G))=\dim(R/I_3(G)).\] 
    Hence, $I_3(G)$ is Cohen-Macaulay. The other implication is well-known.
\end{proof}

\section{Alexander dual of the \texorpdfstring{$3$-path ideals of trees}{}}\label{section 4}

In this section, we study the Alexander dual ideal of the $3$-path ideal of a tree. Although the notion of Alexander duality comes from the theory of simplicial complexes, we consider the equivalent algebraic definition in terms of Stanley-Reisner ideals. As there is a one-to-one correspondence between the Stanley-Reisner ideals and square-free monomial ideals, one can define the Alexander dual of a square-free monomial ideal in the following way.

\begin{definition}\label{alexander}{\rm
    Let $I=\langle x_{11}\cdots x_{1l_1},\ldots,x_{q1}\cdots x_{ql_q}\rangle$ be a square-free monomial ideal in a polynomial ring $R$. The {\it Alexander dual} (in short, \textit{dual}) $I^{\vee}$ of $I$ is defined as
\[
 I^{\vee}:= \langle x_{11},\ldots, x_{1l_1}\rangle\cap\dots\cap\langle x_{q1},\ldots ,x_{ql_q}\rangle.
\]
}
\end{definition}

Recall that, for a hypergraph $\H$, the minimal vertex covers are those vertex covers of $\H$ that are minimal with respect to inclusion. The collection of all minimal primes of $I(\H)$ are in one-to-one correspondence with the minimal vertex covers of $\H$. Therefore,
 \[I(\H)^{\vee}=\langle x_{i_1}\cdots x_{i_t}\mid \{x_{i_1},\ldots,x_{i_t}\} \text{ is a minimal vertex cover of } \H \rangle.\]
 Now, given a hypergraph $\H$, let us define a new hypergraph $\H^{\vee}$, where $V(\H^{\vee})=V(\H)$ and $E(\H^{\vee})=\{\E\subseteq V(\H)\mid \E\text{ is a minimal vertex cover of }\H\}$. Then $I(\H)^{\vee}=I(\H^{\vee})$. We shall make use of this identification of the dual of a square-free monomial ideal throughout this section.

Recently, Moradi and Khosh-Ahang \cite{MKA} defined the notion of vertex splittable monomial ideals, which 
provides an algebraic analogue of the vertex decomposable property of an abstract simplicial complex. In particular, a simplicial complex is vertex decomposable if and only if the dual of its Stanley-Reisner ideal is vertex splittable (\cite[Theorem 2.3]{MKA}). Our aim here is to show that when $G$ is a tree, the Alexander dual ideal $I(\P_3(G))^{\vee}$ is vertex splittable.
    
\begin{definition}\label{vsp defn}
\normalfont
 Let $I$ be a monomial ideal in the polynomial ring $R=\mathbb K[x_1,\ldots,x_n]$. We say that $I$ is {\it vertex splittable} if $I$ can be obtained by the following recursive procedure.
 \begin{enumerate}
  \item[(i)] If $I=\langle m\rangle$, where $m$ is a monomial or $I=\l 0\r$ or $I=R$, then $I$ is a vertex splittable ideal.
  
  \item[(ii)] If there exists a variable $x_i$ and two vertex splittable ideals $I_1$ and $I_2$ of $\mathbb K[x_1,\ldots,\widehat{x_i},\ldots,x_n]$ such that $I=x_iI_1+I_2$ with $I_2\subseteq I_1$ and the minimal generators of $I$ is the disjoint union of the minimal generators of $x_iI_1$ and $I_2$, then $I$ is a vertex splittable ideal. 
 \end{enumerate}
\end{definition}

\begin{remark}
    If $I$ is a vertex splittable ideal, then the vertex $x_i$ in \Cref{vsp defn} is called the splitting vertex of $I$.
\end{remark}

We frequently use the following well-known formula due to N. Terai \cite[Theorem 2.1]{Terai}, which relates the regularity and projective dimension of a square-free monomial ideal and its dual.

\begin{theorem}[Terai's formula]
Let $I$ be a square-free monomial ideal in $R$. Then \[\pd(R/I)=\reg(I^{\vee}).\]
\end{theorem}

\noindent
We now recall a few preliminary results related to the vertex splittable property of monomial ideals that we shall use in this section.

\begin{proposition}\label{two product splittable}
    Let $I\subseteq R_1=\K[x_1,\ldots,x_n]$ and $J\subseteq R_2=\K[y_1,\ldots,y_m]$ be two vertex splittable monomial ideals. Then $IJ\subseteq R_1\otimes_{\K} R_2$ is also a vertex splittable ideal.
\end{proposition}
\begin{proof}
    If one of $I$ and $J$ is $\l0\r$ or $\l1\r$, then we are done. Also, if both $I$ and $J$ are principal monomial ideals, then $IJ$ is also a principal monomial ideal, and so is vertex splittable. Now, without loss of generality, we may assume that $I$ is not a principal ideal. We proceed by induction on $n+m$. The cases $n+m\leq 2$ easily follow from the previous arguments. So, we assume $n+m> 2$. Since $I$ is a vertex splittable ideal and not generated by a single monomial, we can write $I=x_iI_1+I_2$, where $I_1,I_2\subseteq \K[x_1,\ldots,\widehat{x_i},\ldots,x_n]$ are vertex splittable ideals, $I_2\subseteq I_1$, and $\G(I)=\G(x_iI_1)\sqcup\G(I_2)$. Now $IJ=x_iI_1J+I_2J$, where $I_2J\subseteq I_1J$ and $\G(IJ)=\G(x_iI_1J)\sqcup\G(I_2J)$. Using the induction hypothesis, we conclude that $I_1J$ and $I_2J$ are vertex splittable ideals. Since no minimal generator of $I_1J$ and $I_2J$ is divisible by $x_i$, $IJ$ is vertex splittable.
\end{proof}

\begin{corollary}\label{forest corollary}
    Let $I_i$ be monomial ideals in the polynomial rings $R_i=\K [x_{i1},\ldots,x_{ir_i}]$ respectively, where $i\in [t]$. If each $I_i$ is a vertex splittable ideal in $R_i$ for $i\in[t]$, then $I=\prod_{i=1}^tI_i$ is a vertex splittable ideal in $R_1\otimes_{\K}\dots\otimes_{\K}R_t$.
\end{corollary}

The following proposition shows that the vertex splittable property of the edge ideal of a hypergraph is preserved even if we delete a vertex from the hypergraph. More generally, we show the following.

\begin{proposition}\label{vertex omit}
    Let $\H$ be a hypergraph such that $I(\H)$ is a vertex splittable ideal. Then for any induced subhypergraph $\H'$ of $\H$, the ideal $I(\H')$ is also vertex splittable.
\end{proposition}
\begin{proof}
    It is enough to show that $I(\H\setminus x)$ is vertex splittable for any $x\in V(\H)$. We proceed by induction on $|V(\H)|$. If $|V(\H)|=1$, then $I(\H\setminus x)$ is the zero ideal, and hence, vertex splittable. Therefore, we may assume that $|V(\H)|\ge 2$.\par 
    
    Let $y\in V(\H)$ be the splitting vertex of $I(\H)$. Consider the two hypergraphs $\H_1$ and $\H_2$, where $V(\H_1)=V(\H_2)=V(\H)\setminus\{y\}$, $E(\H_1)=\{\E\setminus\{y\}\mid \E\in E(\H)\text{ with }y\in\E\}$, and $E(\H_2)=\{\E\in E(\H)\mid y\notin \E\}$. Note that $I(\H)=y I(\H_1)+I(\H_2)$, where $I(\H_2)\subseteq I(\H_1)$ and both $I(\H_1)$ and $I(\H_2)$ are vertex splittable ideals. If $x=y$, then $I(\H\setminus x)=I(\H_2)$, and hence, $I(\H\setminus x)$ is a vertex splittable ideal. If $x\neq y$, then $I(\H\setminus x)=y I(\H_1\setminus x)+I(\H_2\setminus x)$. It is easy to see that $I(\H_2\setminus x)\subseteq I(\H_1\setminus x)$. Also, by the induction hypothesis, both $I(\H_1\setminus x)$ and $I(\H_2\setminus x)$ are vertex splittable ideals.  Therefore, $I(\H\setminus x)$ is a vertex splittable ideal.
\end{proof}

The following lemma and its applications are essential in our proof of \Cref{main theorem1} below.

\begin{lemma}\label{dual vill lemma}\textup{(cf. \cite[Lemma 2.6]{MMVV})}
    Let $\H$ be a hypergraph and let $f=\prod_{j=1}^rx_{i_j}$ be a square-free monomial in $R$. Then 
    \[
    (I(\H)^{\vee}:f)^{\vee}=I(\H\setminus \{x_{i_1},\ldots,x_{i_r}\}).
    \]
    In particular, if $x\in V(\H)$, then $(I(\H)^{\vee}:x)^{\vee}=I(\H\setminus x)$.
\end{lemma}

\begin{proposition}\label{new lemma}
    Let $\H$ be a hypergraph and $x_{i_1},\ldots,x_{i_{r+1}}\in V(\H)$ for some $r\ge 1$. Suppose $\H_1$ and $\H_2$ are two hypergraphs such that $I(\H_1)=(I(\H):x_{i_1}\cdots x_{i_r})$ and $I(\H_2)=I(\H\setminus x_{i_{r+1}})$. Then 
    \[
    I(\H_1\setminus x_{i_{r+1}})^{\vee}=I(\H_2^{\vee}\setminus\{x_{i_1},\ldots,x_{i_r}\}).
    \]
\end{proposition}
\begin{proof}
    It is easy to see that $I(\H_1\setminus x_{i_{r+1}})=(I(\H_2):f)=(I(\H_2^{\vee})^{\vee}:f)$, where $f=\prod_{j=1}^rx_{i_j}$. Hence, by \Cref{dual vill lemma}, $I(\H_1\setminus x_{i_{r+1}})^{\vee}=I(\H_2^{\vee}\setminus\{x_{i_1},\ldots,x_{i_r}\})$.
\end{proof}

Using \Cref{dual vill lemma}, we also obtain the following.

\begin{proposition}\label{splitting variable}
    Let $\H$ be a hypergraph, and $x$ be any vertex of $\H$. Suppose $\H_1$ and $\H_2$ are two hypergraphs with $V(\H_1)=V(\H_2)=V(\H)\setminus \{x\}$, and
    \begin{align*}
        E(\H_1)=\{\E\setminus\{x\}\mid \E\in E(\H^{\vee})\text{ with }x\in\E\},
        E(\H_2)=\{\E\in E(\H^{\vee})\mid x\notin \E\}.
    \end{align*} 
    Then $I(\H)^{\vee}=xI(\H_1)+I(\H_2).$ Moreover, if $I(\H_2)\subseteq I(\H_1)$, then $I(\H_1)=I(\H\setminus x)^{\vee}$ and $I(\H_2)=(I(\H):x)^{\vee}$.
\end{proposition}

\begin{proof}
    Note that $E(\H_1)$ consists of the edges $\E\setminus\{x\}$, where $\E$ is a minimal vertex cover of $\H$ and $x\in\E$. Also, $E(\H_2)$ consists of all the minimal vertex covers of $\H$ which does not contain $x$. Hence, it is easy to observe that $I(\H)^{\vee}=xI(\H_1)+I(\H_2)$.
\par 
    Now if $I(\H_2)\subseteq I(\H_1)$, then $I(\H_1)=(I(\H)^{\vee}:x)$. Hence, $I(\H_1)^{\vee}=(I(\H)^{\vee}:x)^{\vee}=I(\H\setminus x)$ (by \Cref{dual vill lemma}). Thus, $I(\H_1)=I(\H\setminus x)^{\vee}$. Also, $I(\H_2)=I(\H^{\vee}\setminus x)=(I(\H):x)^{\vee}$ (again by \Cref{dual vill lemma}).
\end{proof}

    \begin{figure}
        \centering
         \begin{tikzpicture}
    [scale=.55]
\draw [fill] (4,6) circle [radius=0.1];
\draw [fill] (3,4) circle [radius=0.1];
\draw [fill] (5,4) circle [radius=0.1];
\draw [fill] (2,2) circle [radius=0.1];
\draw [fill] (4,2) circle [radius=0.1];
\draw [fill] (6,2) circle [radius=0.1];
\draw [fill] (3,0) circle [radius=0.1];
\draw [fill] (5,0) circle [radius=0.1];
\node at (4.8,6) {$x$};
\node at (2,4) {$z_1$};
\node at (6,4) {$z_2$};
\node at (1,2) {$z_3$};
\node at (4.5,2) {$z_4$};
\node at (6.5,2) {$z_5$};
\node at (2,0) {$z_6$};
\node at (5.8,0) {$z_7$};
\draw (4,6)--(3,4)--(4,2)--(3,0);
\draw (3,4)--(2,2);
\draw (4,2)--(5,0);
\draw (4,6)--(5,4)--(6,2);
\end{tikzpicture}
        \caption{A rooted tree $T$.}
        \label{fig:enter-label1}
    \end{figure}
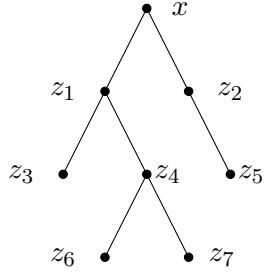

For a graph $G$, note that $I_3(G)=I(\H)$, where $\H=\P_3(G)$, and hence, $I_3(G)^{\vee}=I(\H)^{\vee}=I(\H^{\vee})$. From now on, we shall write $T$ instead of $G$ to denote a tree. Recall that a \textit{tree} $T$ is a finite simple graph such that there exists a unique path between
any two distinct vertices. Fix any vertex $x\in V(T)$. Then $T$ can be realised as a rooted tree with the vertex $x$ treated as a root. For any $z\in V(T)$, if $x=w_{i_1},w_{i_2},\ldots,w_{i_{r-1}},w_{i_r}=z$ is the unique path between $x$ and $z$, then we say that $z$ has {\it level} $r$ and denote it by $\mathrm{level}(z)=r$. Moreover, we define the height of $T$ to be
\[
\mathrm{height}(T)=\max_{z\in V(T)}\mathrm{level}(z).
\]
\noindent
As an example, consider the rooted tree $T$ in \Cref{fig:enter-label1} with the vertex $x$ as a root. Then we have $\mathrm{level}(z_i)=1$ for $i=1,2$, $\mathrm{level}(z_i)=2$ for $i=3,4,5$, and $\mathrm{level}(z_i)=3$ for $i=6,7$. Also, $\mathrm{height}(T)=3$.

Now we proceed to prove the main theorem of this section. 
\begin{theorem}\label{main theorem1}
    Let $T$ be a tree. 
    Then the Alexander dual ideal $I_3(T)^{\vee}$ is vertex splittable.
\end{theorem}
\begin{proof}
       We prove this by induction on $|V(T)|$. If $|V(T)|\le 3$, then $I_3(T)$ is either a zero ideal or generated by a monomial. Then, clearly $I_3(T)^{\vee}$ is a vertex splittable ideal. Therefore, we may assume that $|V(T)|\ge 4$. First we consider the case when $T$ is a star graph with $V(T)=\{w,v_1,\ldots, v_m\}$ and $E(T)=\{\{w,v_i\}\mid i\in [m]\}$. Then $I_3(T)^{\vee}=wJ_1+J_2$, where $J_1=R$ and $J_2$ is the vertex cover ideal of the complete graph on the vertex set $\{v_1,\ldots,v_m\}$. By \cite[Corollary 3.8]{MKA}, $J_2$ is a vertex splittable ideal. Also, $J_2\subseteq J_1$ and by definition, $J_1$ is a vertex splittable ideal. Hence, $I_3(T)^{\vee}$ is a vertex splittable ideal. Now, let us assume that $T$ is not a star graph. 
\begin{figure}
    \centering
     \begin{tikzpicture}
    [scale=.45]
        \draw [fill] (4.5,9) circle [radius=0.1]; 
        \draw [fill] (3,6) circle [radius=0.1]; 
        \draw [fill] (2,4) circle [radius=0.1]; 
        \draw [fill] (1,2) circle [radius=0.1]; 
        \draw [fill] (2.5,2) circle [radius=0.1]; 
        \draw [fill] (5,2) circle [radius=0.1]; 
        \draw [fill] (0,0) circle [radius=0.1]; 
        \draw [fill] (1.3,0) circle [radius=0.1]; 
        \draw [fill] (4.5,0) circle [radius=0.1]; 
        \node at (4.5,9.5) {$w$};
        \node at (2.5,6) {$z_1$};
        \node at (1.5,4) {$b$};
        \node at (0.5,2) {$y$};
        \node at (2.5,1.5) {$z_2$};
        \node at (5,1.5) {$z_r$};
        \node at (0,-0.5) {$x_1$};
        \node at (1.3,-0.5) {$x_2$};
        \node at (4.5,-0.5) {$x_l$};
        \node at (2.9,0) {$\cdots$};
        \node at (3.75,2) {$\cdots$};
        \draw (3,6)--(2,4)--(1,2)--(0,0);
        \draw (1,2)--(1.3,0);
        \draw (1,2)--(4.5,0);
        \draw (2,4)--(2.5,2);
        \draw (2,4)--(5,2);
        \draw[] (4.5,9)--(4,8);
        \draw[] (4.5,9)--(4.5,8);
        \draw[] (4.5,9)--(5,8);
        \draw[dashed] (4.5,8)--(4.5,7);
        \draw[dashed] (5,8)--(5.5,7);
        \draw[dashed] (4,8)--(3,6);
    \end{tikzpicture}
   
    \caption{A tree with splitting vertices $y$ and $b$.}
    \label{fig:enter-label}
\end{figure}
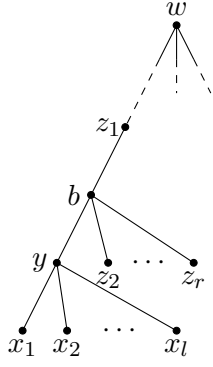

    First, consider the tree $T$ as a rooted tree with some vertex, say $w$, as a root. Choose any vertex $y$ in $T$ such that $\level(y)=\mathrm{height}(T)-1$. Then $N_T(y)=\{x_1,\ldots,x_l,b\}$, where $x_1,\ldots,x_l$ are leaves of $T$ with $l\ge 1$. Since $T$ is not a star graph, $b$ is not a leaf. We have the following two cases:
    
\noindent 
    {\bf Case I}: Assume that $l\ge 2$. In this case, by \Cref{splitting variable}, $I(\H)^{\vee}=yI(\H_1)+I(\H_2)$, where $\H_1$ and $\H_2$ are two hypergraphs both on the vertex set $V(T)\setminus\{y\}$ with edge sets $E(\H_1)=\{\E\setminus\{y\}\mid \E\in E(\H^{\vee})\text{ with }y\in\E\}$ and $E(\H_2)=\{\E\in E(\H^{\vee})\mid y\notin \E\}$. Note that $\G(I(\H)^{\vee})=\G(yI(\H_1))\sqcup \G(I(\H_2))$.

    \noindent
    {\bf Claim 1}: $I(\H_2)\subseteq I(\H_1)$.

    \noindent{\bf Proof of Claim 1}: Let $\E\in E(\H^{\vee})$ be such that $y\notin \E$. Observe that if $x_1,x_2\notin\E$, then $\E\cap\{x_1,y,x_2\}=\emptyset$, a contradiction since $\{x_1,y,x_2\}\in E(\H)$. Without any loss of generality, let $x_1\in \E$. We show that $\E'=(\E\setminus\{x_1\})\cup\{y\}$ is a vertex cover of $\H$. Indeed, if $\mathcal F\in E(\H)$ such that $x_1\notin\F$, then $\E\cap\F\subseteq \E'\cap\F$. Hence, $\E'\cap\F\neq\emptyset$ since $\E\cap\F\neq\emptyset$.  Also if $x_1\in\F$, then $y\in\E'\cap\F$, since $x_1\in\F$ implies $y\in \F$. Thus, $\E'$ is a vertex cover of $\H$. Now if $\E''\in E(\H^{\vee})$ such that $\E''\subseteq \E'$, then $y\in\E''$. Because otherwise, $\E''\subsetneq\E$ is a vertex cover of $\H$, a contradiction. Therefore, $\x_{\E''\setminus\{y\}}\in I(\H_1)$, and $\x_{\E''\setminus\{y\}}$ divides $\x_{\E}$, where $\x_{\E}\in I(\H_2)$ is arbitrarily chosen. Consequently, $I(\H_2)\subseteq I(\H_1)$.

    Now, by \Cref{splitting variable} we have $I(\H_1)=I(\H\setminus y)^{\vee}$ and $I(\H_2)=(I(\H):y)^{\vee}$. Note that the ideal $I(\H\setminus y)$ is the $3$-path ideal of the tree $T\setminus\{y,x_1,\ldots,x_l\}$. Therefore, by the induction hypothesis, $I(\H_1)$ is vertex splittable. Now, it remains to show that the ideal $(I(\H):y)^{\vee}$ is vertex splittable. So, we consider the ideal $(I(\H):y)$. Let $\H_3$ be the hypergraph on the vertex set $V(T)\setminus \{y\}$ such that $I(\H_3)=(I(\H):y)$. Let $N_T(b)=\{y,z_1,\ldots,z_r\}$. Then 
    \[
    I(\H_3)=(I(\H):y)=\l x_ix_j\mid 1\le i<j\le l \r+\l b x_i,bz_j\mid i\in[l],j\in[r]\r+I_3(T\setminus N_T[y]).
    \]
     Again by \Cref{splitting variable}, $I(\H_3)^{\vee}=bI(\H_4)+I(\H_5)$ where $\H_4$ and $\H_5$ are two hypergraphs both on the vertex set $V(T)\setminus \{y,b\}$ with edge sets $E(\H_4)=\{\E\setminus\{b\}\mid \E\in E(\H_3^{\vee})\text{ with }b\in\E\}$, and $E(\H_5)=\{ \E\in E(\H_3^{\vee})\mid b\notin\E\}$. Note that $\G(I(\H_3)^{\vee})=\G(bI(\H_4))\sqcup \G(I(\H_5))$.

    \noindent {\bf Claim 2}: $I(\H_5)\subseteq I(\H_4)$.
    
    \noindent  {\bf Proof of Claim 2}: Let $\E\in E(\H_3^{\vee})$ be such that $b\notin\E$. In that case $x_i\in\E$ for each $i\in [l]$. We show that $\E'=(\E\setminus\{x_1\})\cup\{b\}$ is a vertex cover of $\H_3$. Indeed, $\E'\cap\{x_i,x_j\}\neq\emptyset$ for $i,j\in[l]$ since either $i\neq 1$ or $j\neq 1$. Also, $b\in\E'\cap \{b,x_i\}$ and $b\in \E'\cap\{b,z_j\}$. Moreover, if $\F\in E(\H_3)$ such that $\F\subseteq V(T)\setminus N_T[y]$, then $\E'\cap\F=\E\cap\F\neq\emptyset$. Thus, $\E'$ is a vertex cover of $\H_3$. Now, if $\E''\in E(\H_3^{\vee})$ such that $\E''\subseteq \E'$, then $b\in \E''$. Because otherwise, $\E''\subsetneq \E$ is a vertex cover of $\H_3$, a contradiction. Therefore, $\x_{\E''\setminus\{b\}}\in I(\H_4)$ and $\x_{\E''\setminus\{b\}}$ divides $\x_{\E}$, where $\x_{\E}\in I(\H_5)$ is arbitrarily chosen. Consequently, $I(\H_5)\subseteq I(\H_4)$. \par 
    
    Now, by \Cref{splitting variable}, we have $I(\H_4)=I(\H_3\setminus b)^{\vee}$ and $I(\H_5)=(I(\H_3):b)^{\vee}$. Note that $(I(\H_3):b)=(I(\H):by)$ and therefore, $(I(\H_3):b)=\l x_1,\ldots,x_l,z_1,\ldots,z_r\r+I_3(T\setminus N_T[y,b])$. Thus, $I(\H_5)=\left(\prod_{i=1}^lx_i\right)\cdot\left(\prod_{j=1}^rz_j\right)\cdot I_3(T\setminus N_T[y,b])^{\vee}$. Observe that $T\setminus N_T[y,b]$ is a forest, and each connected component is a tree with lesser number of vertices than that of $T$. Therefore, by the induction hypothesis, and by \Cref{forest corollary}, $I_3(T\setminus N_T[y,b])^{\vee}$ is a vertex splittable ideal, and hence, $I(\H_5)$ is also a vertex splittable ideal. Now, $I(\H_3\setminus b)=\l x_ix_j\mid 1\le i<j\le l \r+I_3(T\setminus N_T[y])$. Note that $J=\l x_ix_j\mid 1\le i<j\le l \r$ is the edge ideal of the complete graph on $l$ vertices. Therefore, by \cite[Corollary 3.8]{MKA}, $J^{\vee}$ is a vertex splittable ideal. Moreover, $T\setminus N_T[y]$ is a forest, and thus proceeding as before, we get that $I_3(T\setminus N_T[y])^{\vee}$ is vertex splittable. Consequently, by \Cref{two product splittable}, $I(\H_4)=J^{\vee}\cdot I_3(T\setminus N_T[y])^{\vee}$ is a vertex splittable ideal. Since $I(\H_3)^{\vee}=bI(\H_4)+I(\H_5)$, we obtain that, $I(\H_3)^{\vee}$ is a vertex splittable ideal. Note that $I(\H_3)^{\vee}=I(\H_2)$, and hence, $I(\H)^{\vee}=yI(\H_1)+I(\H_2)$ is a vertex splittable ideal.

\noindent
    {\bf Case II}: Assume that $l=1$.
    We shall proceed similarly as in Case I, but the choice of splitting vertices will be slightly different.
    Let $N_T(b)=\{y,z_1,\ldots,z_r\}$. Then by \Cref{splitting variable}, $I(\H)^{\vee}=bI(\H_6)+I(\H_7)$, where $\H_6$ and $\H_7$ are two hypergraphs both on the vertex set $V(T)\setminus\{b\}$ with edge sets $E(\H_6)=\{\E\setminus\{b\}\mid \E\in E(\H^{\vee})\text{ with }b\in\E\}$ and $E(\H_7)=\{\E\in E(\H^{\vee})\mid b\notin \E\}$. Now $\G(I(\H)^{\vee})=\G(bI(\H_6))\sqcup \G(I(\H_7))$. Also, proceeding as in Case I we have $I(\H_7)\subseteq I(\H_6)$. Therefore, by \Cref{splitting variable}, 
    \[
    I(\H_6)=I(\H\setminus b)^{\vee}\text{ and }I(\H_7)=(I(\H):b)^{\vee}.
    \]
    The ideal $I(\H\setminus b)$ is the $3$-path ideal of the graph $T\setminus b$. Observe that $T\setminus b$ is a forest, and thus proceeding as before, we see that $I(\H_6)=I(\H\setminus b)^{\vee}$ is a vertex splittable ideal. Now we consider the ideal $(I(\H):b)$. Let $N_T(z_i)=\{w_{i1},\ldots,w_{ip_i}\}$ for each $i\in[r]$. Let $\H_8$ be the hypergraph on the vertex set $V(T)\setminus \{b\}$ such that $I(\H_8)=(I(\H):b)$. Then
    \begin{align*}
        I(\H_8)=(I(\H):b)=&\l yx_1,yz_1,\ldots,yz_r\r+\l z_iz_j\mid 1\le i<j\le r \r\\
        &+\l z_iw_{ij}\mid i\in[r],j\in [p_i] \r+I_3(T\setminus N_T[y,b]).
    \end{align*}
     Now by \Cref{splitting variable}, $I(\H_8)^{\vee}=yI(\H_9)+I(\H_{10})$, where $\H_9$ and $\H_{10}$ are two hypergraphs both on the vertex set $V(T)\setminus \{y,b\}$ with edge sets $E(\H_9)=\{\E\setminus\{y\}\mid \E\in E(\H_8^{\vee})\text{ with }y\in\E\}$, and $E(\H_{10})=\{\E\in E(\H_8^{\vee})\mid y\notin\E\}$. It is easy to see that $\G(I(\H_8)^{\vee})=\G(yI(\H_9))\sqcup \G(I(\H_{10}))$. Also, as before $I(\H_{10})\subseteq I(\H_9)$, and thus, again by \Cref{splitting variable},
     \[
     I(\H_9)=I(\H_8\setminus y)^{\vee} \text{ and } I(\H_{10})=(I(\H_8):y)^{\vee}.
     \]

\noindent
    Now $(I(\H_8):y)=(I(\H):by)$, and hence, $(I(\H_8):y)=\l x_1,z_1,\ldots,z_r \r+I_3(T\setminus N_T[y,b])$. Therefore, $I(\H_{10})=x_1\cdot \prod_{i=1}^rz_i\cdot I_3(T\setminus N_T[y,b])^{\vee}$. As in Case I, $T\setminus N_{T}[y,b]$ is a forest, and thus, $I_3(T\setminus N_T[y,b])^{\vee}$ is a vertex splittable ideal. Hence, by \Cref{forest corollary} $I(\H_{10})$ is a vertex splittable ideal.
    Now, using \Cref{new lemma} we obtain $I(\H_8\setminus y)^{\vee}=I((\H\setminus y)^{\vee}\setminus b)$. Note that $I(\H\setminus y)$ is the $3$-path ideal of the tree $T\setminus \{y,x_1\}$. Therefore, by the induction hypothesis $I(\H\setminus y)^{\vee}$ is a vertex splittable ideal. Then using \Cref{vertex omit}, we see that $I((\H\setminus y)^{\vee}\setminus b)$ is a vertex splittable ideal. Consequently, $I(\H_9)$ is a vertex splittable ideal and this implies $I(\H_8)^{\vee}=yI(\H_9)+I(\H_{10})$ is a vertex splittable ideal. Moreover, since $I(\H_8)^{\vee}=I(\H_7)$ and $I(\H)^{\vee}=bI(\H_6)+I(\H_7)$, we conclude that $I(\H)^{\vee}$ is also a vertex splittable ideal. This completes the proof.  
\end{proof}

\begin{remark}\normalfont
    The vertex splittability of an ideal is a stronger property, which implies that the ideal is componentwise linear. Later, we will observe (in \cref{Example1}) that there are chordal graphs $G$ for which the dual of $I_3(G)$ is not even componentwise linear. Moreover, in the case of trees, it is not possible to generalize this result for higher $t$-path ideals (see \Cref{caterpillar not vs}).
\end{remark}

Recently, a generalization of the independence complex of a graph and the properties of its Stanley-Reisner ideal are discussed in \cite{ADGRS}. Recall that, for a positive integer $r$, a subset $W$ of the vertex set of a graph $G$ is called {\it $r$-independent} if each connected component of the induced subgraph $G[W]$ has vertex cardinality at most $r$. The collection of all $r$-independent sets form the {\it $r$-independence complex} of $G$, and it is denoted by $\mathrm{Ind}_r(G)$. Note that the Stanley-Reisner ideal of $\mathrm{Ind}_r(G)$ is the ideal $J_{r+1}(G)=\left\l \prod_{j=1}^{r+1}x_{i_j}\mid G[\{x_{i_1},\ldots,x_{i_{r+1}}\}]\text{ is connected }\right\r$. Abdelmalek et. al. conjectured in \cite[Conjecture 3.15]{ADGRS} that if $T$ is a tree, then the complex $\mathrm{Ind}_r(T)$ is vertex decomposable for any positive integer $r$. We prove this conjecture for $r=2$.

\begin{corollary}\label{vdtree}
    If $T$ is a tree, then $\mathrm{Ind}_2(T)$ is vertex decomposable.
\end{corollary}
\begin{proof}
    Note that $\mathrm{Ind}_2(T)$ is the Stanley-Reisner complex of the ideal $I(\P_3(T))=I_3(T)$. Since $I_3(T)^{\vee}$ is vertex splittable, by \cite[Theorem 2.3]{MKA}, $\mathrm{Ind}_2(T)$ is a vertex decomposable simplicial complex.   
\end{proof}
 
\begin{remark}\normalfont
    For square-free monomial ideals, the vertex splittable property is equivalent to the vertex decomposability of the corresponding simplicial complexes. Notably, the notion of vertex decomposability of simplicial complexes is independent of the choice of the base field. On the other hand, the (sequentially) Cohen-Macaulayness of a simplicial complex depends on the choice of the base field. The notion of (sequentially) Cohen-Macaulay property of a simplicial complex is equivalent to the (sequentially) Cohen-Macaulayness of the corresponding Stanley-Reisner ring. Moreover, vertex decomposable simplicial complexes are sequentially Cohen-Macaulay.  Therefore, it follows from \Cref{main theorem1} that when $G$ is a tree, the (sequentially) Cohen-Macaulay property of $I_3(G)$ is independent of the choice of the base field. 
\end{remark}

\section{\texorpdfstring{$t$-path}{} ideals of caterpillar graphs and some concluding remarks}\label{section 6}

In this section, we try to extend the results that we have obtained in the previous sections to $t$-path ideals for all $t\geq 3$. Regularity, projective dimension, and Betti numbers of all $t$-path ideals are only known for paths and cycles (see \cite{AlilooeeFaridi2015, AlilooeeFaridi2018}). Here, we give an analogous formula for the regularity of all $t$-path ideals of caterpillar graphs (a special class of trees) when compared to the $3$-path ideals of chordal graphs. Moreover, we will show that our main results can not be generalized to higher $t$-path ideals, even in case of trees. At the end, we give some concluding remarks and questions.

A \textit{caterpillar} graph is a tree from which, if we remove all the pendant vertices, then what remains is a path graph. For our convenience, we write the vertex set and the edge set of a caterpillar graph as follows:
\begin{align*}
V(G)&=\{x_1,x_n,x_i,y_{i,j}\mid 2\le i\le n-1,j\in [r_i]\},\\
E(G)&=\{\{x_1,x_2\},\{x_i,x_{i+1}\},\{x_i,y_{i,j}\}\mid 2\le i\le n-1,j\in [r_i]\}.
\end{align*}
Note that for some $i\in\{2,\ldots,n-1\}$, $r_i$ may be zero. We now extend \Cref{result1} for all $t$-path ideals of a caterpillar graph.
\begin{theorem}\label{caterpillar reg}
    Let $G$ be a caterpillar graph. Then for all $t\ge 2$, $\reg(R/I_t(G))= (t-1)\nu_t(G)$, where $\nu_t(G)$ denotes the $t$-path induced matching of $G$.
\end{theorem}

\begin{proof}
Let $G$ be a caterpillar graph with $V(G)$ and $E(G)$ given as above. First we prove the inequality $\reg(R/I_t(G))\le (t-1)\nu_t(G)$ by induction on $|V(G)|$. If $|V(G)|\le t-1$, then $I_t(G)=\l 0\r$ and hence, $\reg(R/I_t(G))=(t-1)\nu_t(G)$. Therefore, we may assume that $|V(G)|\ge t$. Now if $n\le t-1$, then again $I_t(G)=\l 0\r$, and hence, $\reg(R/I_t(G))=(t-1)\nu_t(G)$. Thus, we may assume that $|V(G)|\ge n\ge t$. In that case it is easy to see that $\nu_t(G)\ge 1$. Now let $J=\l m\in I_t(G)\mid x_1\mid m \r$, and $K=I_t(G\setminus x_1)$. Then $J+K=I_t(G)$. Note that $J=\left(\prod_{i=1}^{t-1}x_i\right)\cdot \left\l\{ x_t,y_{t-1,1},\ldots,y_{t-1,r_{t-1}} \}\right\r$. Therefore, we can write 
\[
J\cap K=\left(\prod_{i=1}^{t-1}x_i\right)\cdot L,
\]
where $L$ is generated by the monomials $\frac{\lcm(m,m')}{\prod_{i=1}^{t-1}x_i}$ such that $m\in J$ and $m'\in K$.

\noindent
{\bf Claim 1:} $L=\l y_{2,j}y_{t-1,k}\mid  j\in[r_2],k\in[r_{t-1}]\r +\l y_{2,j}x_t\mid j\in[r_2]\r +\l x_ty_{t,l}\mid l\in[r_{t}]\r \newline 
    \text{\hspace{5.8cm}}+ \l x_tx_{t+1}\r +\sum_{w\in S}wI_t(G\setminus(N_G[x_1,\ldots,x_t]\setminus\{x_{t+1}\}))$,

\noindent
where $S=\{x_t,y_{t-1,1},\ldots,y_{t-1,r_{t-1}}\}$.

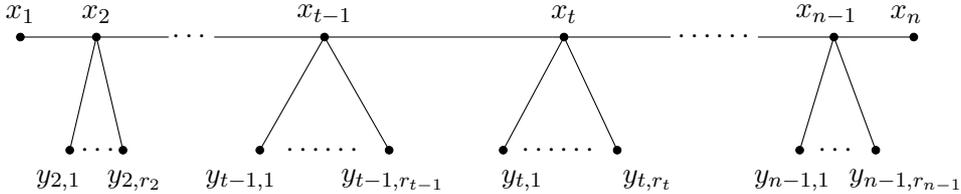
\begin{figure}[!ht]
\centering
\begin{tikzpicture}
[scale=.50]
\draw [fill] (-2,3) circle [radius=0.1];
\draw [fill] (0,3) circle [radius=0.1];
\draw [fill] (6,3) circle [radius=0.1];
\draw [fill] (12.3,3) circle [radius=0.1];
\draw [fill] (19.4,3) circle [radius=0.1];
\draw [fill] (18.5,0) circle [radius=0.1];
\draw [fill] (20.5,0) circle [radius=0.1];
\draw [fill] (21.5,3) circle [radius=0.1];
\draw [fill] (-0.7,0) circle [radius=0.1];
\draw [fill] (0.7,0) circle [radius=0.1];
\draw [fill] (4.3,0) circle [radius=0.1];
\draw [fill] (7.7,0) circle [radius=0.1];
\draw [fill] (10.7,0) circle [radius=0.1];
\draw [fill] (13.7,0) circle [radius=0.1];
\node at (-2,3.6) {$x_1$};
\node at (0,3.6) {$x_2$};
\node at (6,3.6) {$x_{t-1}$};
\node at (19.2,3.6) {$x_{n-1}$};
\node at (21.3,3.6) {$x_{n}$};
\node at (12.3,3.6) {$x_{t}$};
\node at (-1,-0.8) {$y_{2,1}$};
\node at (1,-0.8) {$y_{2,r_2}$};
\node at (3.8,-0.8) {$y_{t-1,1}$};
\node at (7.8,-0.8) {$y_{t-1,r_{t-1}}$};
\node at (11.2,-0.8) {$y_{t,1}$};
\node at (14.5,-0.8) {$y_{t,r_{t}}$};
\node at (18.3,-0.8) {$y_{n-1,1}$};
\node at (21.3,-0.8) {$y_{n-1,r_{n-1}}$};
 \node at (0.1,0) {$\cdots$};
 \node at (5.5,0) {$\cdots$};
 \node at (6.5,0) {$\cdots$};
 \node at (11.7,0) {$\cdots$};
 \node at (12.7,0) {$\cdots$};
 \node at (15.8,3) {$\cdots$};
 \node at (16.8,3) {$\cdots$};
 \node at (19.5,0) {$\cdots$};
 \node at (2.5,3) {$\cdots$};
 \draw (-2,3)--(0,3)--(-0.7,0);
 \draw (0,3)--(0.7,0);
\draw (6,3)--(4.3,0);
\draw (6,3)--(7.7,0);
\draw (12.3,3)--(10.7,0);
\draw (12.3,3)--(13.7,0);
\draw (0,3)--(1.9,3);
\draw (3.1,3)--(6,3);
\draw (6,3)--(12,3);
 \draw (12,3)--(15.1,3);
 \draw (18.5,0)--(19.4,3)--(21.5,3);
 \draw (17.4,3)--(19.4,3)--(20.5,0);
\end{tikzpicture}\caption{A caterpillar graph $G$.}\label{figure 122465}
\end{figure}

\noindent
{\bf Proof of Claim 1:}  Let $L'$ denote the monomial ideal in the right hand side of the above expression. We proceed to show that $L'\subseteq L$. Note that if $w\in S$ and $m\in I_t(G\setminus(N_G[x_1,\ldots,x_t]\setminus\{x_{t+1}\}))$, then we can write $\left(\prod_{i=1}^{t-1}x_i\right)\cdot wm=\lcm\left(\prod_{i=1}^{t-1}x_i\cdot w,m\right)$, where $\left(\prod_{i=1}^{t-1}x_i\right)\cdot w\in J$ and $m\in K$. So we get $wm\in L$. Now,
for $j\in [r_2]$ and $k\in [r_{t-1}]$, $\left(\prod_{i=1}^{t-1}x_i\right)\cdot y_{2,j}y_{t-1,k}=\lcm\left(\prod_{i=1}^{t-1}x_i\cdot y_{t-1,k},y_{2,j}y_{t-1,k}\cdot \prod_{i=2}^{t-1}x_i\right)$, where $\left(\prod_{i=1}^{t-1}x_i\right)\cdot y_{t-1,k}\in J$ and $y_{2,j}y_{t-1,k}\cdot \left(\prod_{i=2}^{t-1}x_i\right)\in K$. Thus, $y_{2,j}y_{t-1,k}\in L$ for $j\in[r_2]$ and $k\in [r_{t-1}]$. By similar arguments, one can show that $x_tx_{t+1}\in L$, for $j\in [r_2]$, $y_{2,j}x_t\in L$, and for $l\in [r_t]$, $ x_ty_{t,l}\in L$. Thus, $L'\subseteq L$.

To prove $L\subseteq L'$, we consider the monomials $\lcm(m,m')$, where $m\in J$ and $m'\in K$. Note that if $x_2\mid m'$, then either $x_ty_{t,l}\mid m'$ for some $l\in [r_t]$, or $x_tx_{t+1}\mid m'$, or $y_{2,j}x_t\mid m'$ for some $j\in [r_2]$, or $y_{2,jy_{t-1,k}}\mid m'$ for some $j\in [r_2]$ and $k\in [r_{t-1}]$. Hence, for such an $m'\in K$, $\frac{\lcm(m,m')}{\prod_{i=1}^{t-1}x_i}\in L'$ for any $m\in J$. Now suppose $x_2\nmid m'$ and $x_i\mid m'$ for some $3\le i\le t$. Then either $x_{t}x_{t+1}\mid m'$, or $x_ty_{t,l}\mid m'$ for some $l\in [r_t]$. Therefore, for such an $m'\in K$, $\frac{\lcm(m,m')}{\prod_{i=1}^{t-1}x_i}\in L'$ for any $m\in J$. Now consider the case when $x_i\nmid m'$ for $2\le i\le t$. Then $m'\in I_t(G\setminus(N_G[x_1,\ldots,x_t]\setminus\{x_{t+1}\}))$. In that case, for any $m\in J$, we see that $\frac{\lcm(m,m')}{\prod_{i=1}^{t-1}x_i}=wm'$, where $w\in\{x_t,y_{t-1,1},\ldots,y_{t-1,r_{t-1}}\}$. Therefore, $L\subseteq L'$. This proves the claim. 

\noindent
{\bf Claim 2}: $\reg(R/L)\le (t-1)\nu_t(G)-(t-2)$.

\noindent
{\bf Proof of Claim 2}: First consider the case when $x_{t-1}$ do not have any leaf, i.e., $r_{t-1}=0$. In that case it is easy to see from the above description of $L$ that $\reg(R/L)=1+\reg(R/(L:x_t))$, where $(L:x_t)=I_t(G\setminus(N_G[x_1,\ldots,x_{t+1}]\setminus\{x_{t+2}\}))$ since $\{x_1,\ldots,x_t\}$ forms a $t$-path in $G$. Hence, by the induction hypothesis, $\reg(R/(L:x_t))\le (t-1)\nu_t(G\setminus(N_G[x_1,\ldots,x_{t+1}]\setminus\{x_{t+2}\}))$. Note that if $A$ denotes a $t$-path induced matching of $G\setminus(N_G[x_1,\ldots,x_{t+1}]\setminus\{x_{t+2}\})$ such that $|A|=\nu_t(G\setminus(N_G[x_1,\ldots,x_{t+1}]\setminus\{x_{t+2}\}))$, then $A\cup\{x_1,\ldots,x_t\}$ is a $t$-path induced matching of $G$. Hence,  $\nu_t(G\setminus(N_G[x_1,\ldots,x_{t+1}]\setminus\{x_{t+2}\}))\le \nu_t(G)-1$. Consequently, $\reg(R/L)\le (t-1)\nu_t(G)-(t-2)$. \par

Now suppose $r_{t-1}\ge 1$. Consider the ideal $L_1=L +\l x_t,y_{t-1,1},\ldots,y_{t-1,r_{t-1}-1} \r$. Note that $ (L_1:y_{t-1,r_{t-1}}) = (L:y_{t-1,r_{t-1}})+\l x_t,y_{t-1,1},\ldots,y_{t-1,r_{t-1}-1}\r$ (by \Cref{colon comma exchange}). Therefore, using the expression of $L$ in Claim 1, we obtain $(L_1:y_{t-1,r_{t-1}})=\l y_{2,j}\mid j\in [r_2]\r+\l x_t,y_{t-1,1},\ldots,y_{t-1,r_{t-1}-1}\r+I_t(G\setminus(N_G[x_1,\ldots,x_t]\setminus\{x_{t+1}\}))$. Hence, by the induction hypothesis, $\reg(R/(L_1:y_{t-1,r_{t-1}}))\le (t-1)\nu_t(G\setminus(N_G[x_1,\ldots,x_t]\setminus\{x_{t+1}\}))$. Note that $\{x_1,\ldots,x_{t-1},y_{t-1,1}\}$ is a $t$-path in $G$. Thus, proceeding as in the previous paragraph, we get $\nu_t(G\setminus(N_G[x_1,\ldots,x_t]\setminus\{x_{t+1}\}))\le \nu_t(G)-1$. Therefore, $\reg(R/(L_1:y_{t-1,r_{t-1}}))\le (t-1)\nu_t(G)-(t-1)$.
Now $ L_1+\l y_{t-1,r_{t-1}}\r=\l x_t,y_{t-1,1},\ldots,y_{t-1,r_{t-1}-1},y_{t-1,r_{t-1}}\r$. Thus, $\reg(R/(L_1+\l y_{t-1,r_{t-1}}\r)=0\le (t-1)\nu_t(G)-(t-1)$, since $\nu_t(G)\ge 1$. We now use \Cref{regularity lemma} to conclude that $\reg(R/L_1)\le (t-1)\nu_t(G)-(t-2)$.

Next consider the ideal $L_2= L+\l x_t,y_{t-1,1},\ldots,y_{t-1,r_{t-1}-2} \r$. Using \Cref{colon comma exchange}, we have $(L_2:y_{t-1,r_{t-1}-1})= (L:y_{t-1,r_{t-1}-1})+\l x_t,y_{t-1,1},\ldots,y_{t-1,r_{t-1}-2} \r$. As before, the expression of $L$ in Claim 1 gives $(L_2:y_{t-1,r_{t-1}-1})=\l y_{2,j}\mid j\in [r_2]\r+ \l x_t,y_{t-1,1},\ldots,y_{t-1,r_{t-1}-2} \r+I_t(G\setminus(N_G[x_1,\ldots,x_t]\setminus\{x_{t+1}\}))$. Hence, by the induction hypothesis, we get $\reg(R/(L_2:y_{t-1,r_{t-1}-1}))\le (t-1)\nu_t(G\setminus(N_G[x_1,\ldots,x_t]\setminus\{x_{t+1}\}))$. Since $\nu_t(G\setminus(N_G[x_1,\ldots,x_t]\setminus\{x_{t+1}\}))\le \nu_t(G)-1$, we have, $\reg(R/(L_2:y_{t-1,r_{t-1}-1}))\le (t-1)\nu_t(G)-(t-1)$. Notice that $L_1= L_2+\l y_{t-1,r_{t-1}-1}\r$. Thus, using \Cref{regularity lemma} we conclude that $\reg(R/L_2)\le (t-1)\nu_t(G)-(t-2)$. Continuing the above process, we finally obtain $\reg(R/ 
L+\l x_t \r)\le  (t-1)\nu_t(G)-(t-2)$. Now $(L:x_t)=\l x_{t+1},y_{2,j},y_{t,l} \r+I_t(G\setminus(N_G[x_1,\ldots,x_{t+1}]\setminus\{x_{t+2}\}))$, where $j\in[r_2]$ and $l\in [r_t]$. By the induction hypothesis, we have $\reg(R/I_t(G\setminus(N_G[x_1,\ldots,x_{t+1}]\setminus\{x_{t+2}\})))\le (t-1)\nu_t(G\setminus(N_G[x_1,\ldots,x_{t+1}]\setminus\{x_{t+2}\}))\le(t-1)( \nu_t(G)-1)$. Therefore, using \Cref{regularity lemma}, we conclude that $\reg(R/L)\le (t-1)\nu_t(G)-(t-2)$, and this completes the proof of Claim 2. 

Now since $J\cap K=\left(\prod_{i=1}^{t-1}x_i\right)\cdot L$, we have $\reg(R/J\cap K)\le (t-1)\nu_t(G)+1$. Also, $\reg(R/J)=t-1\le (t-1)\nu_t(G)$ since $\nu_t(G)\ge 1$. Moreover, by the induction hypothesis, $\reg(R/K)\le (t-1)\nu_t(G\setminus x_1)\le (t-1)\nu_t(G)$. Hence, using \Cref{regularity lemma1} we obtain, $\reg(R/I_t(G)= \reg(R/J+K)\le (t-1)\nu_t(G)$. Note that by \Cref{lower bound} we have, $\reg(R/I_t(G)\ge (t-1)\nu_t(G)$. Therefore, $\reg(R/I_t(G))= (t-1)\nu_t(G)$ and this completes the proof of the theorem.
\end{proof}

We now discuss about a possible generalization of one of our main results in \Cref{section 3}. More specifically, we proved in \Cref{result1} that $\reg(R/I_3(G))=2\nu_{3}(G)$ for any chordal graph $G$. Also, for $2$-path ideals (equivalently, edge ideals) of graphs, it is well-known that $\reg(R/I_2(G))=\nu_{2}(G)$ when $G$ is chordal. Again, for any caterpillar graph $G$, we have shown in \Cref{caterpillar reg} that the formula $\reg(R/I_t(G))=(t-1)\nu_t(G)$ hold for all $t\geq 2$. Thus, for the class of chordal graphs, one can expect to ask the following:

\begin{question}\label{qreg}
    Let $G$ be any chordal graph. Then, is it true that $\reg(R/I_t(G))=(t-1)\nu_{t}(G)$ for all $t\geq 4$? If not, does the equality hold for trees?
\end{question}

The above question has a negative answer even for the class of trees. In support, we show in the following theorem that $\reg(R/I_t(T))$ can be arbitrarily larger than $(t-1)\nu_{t}(T)$ for any given $t\geq 4$, where $T$ is a tree.

\begin{theorem}\label{reg tree}
     For a tree $T$, the difference between $\reg(R/I_t(T))$ and the quantity $(t-1)\nu_t(T)$ can be arbitrarily large for each $t\geq 4$.
\end{theorem}
\begin{proof}
\begin{figure}[!ht]
\centering
\begin{tikzpicture}
[scale=.45]
    \draw [fill] (-4,5) circle [radius=0.1];
    \draw [fill] (-1,5) circle [radius=0.1];
    \draw [fill] (6,5) circle [radius=0.1];
    \draw [fill] (9,5) circle [radius=0.1];
    \draw [fill] (10.8,6.2) circle [radius=0.1];
    \draw [fill] (12.5,7.3) circle [radius=0.1];
    \draw [fill] (10.8,3.8) circle [radius=0.1];
    \draw [fill] (12.5,2.7) circle [radius=0.1];
\node at (-4,5.8) {$x_1$};
\node at (-1,5.8) {$x_2$};
\node at (6,5.8) {$x_{t-3}$};
\node at (8.6,5.8) {$x_{t-2}$};
\node at (10.4,7) {$w_{11}$};
\node at (13,7.8) {$w_{12}$};
\node at (10.4,3.2) {$w_{k1}$};
\node at (13,2.2) {$w_{k2}$};
\node at (1.9,5) {$\cdots$};
\node at (3.1,5) {$\cdots$};
\node at (12.5,3.9) {$\vdots$};
\node at (12.5,5.1) {$\vdots$};
\node at (12.5,6.3) {$\vdots$};
\node at (10.8,5.2) {$\vdots$};
\draw (-4,5)--(-1,5)--(1,5);
\draw (4,5)--(6,5)--(9,5)--(10.8,6.2)--(12.5,7.3);
\draw (12.5,2.7)--(10.8,3.8)--(9,5);
    \end{tikzpicture}\caption{The graph $T_t'$.}\label{figure T'}
\end{figure}
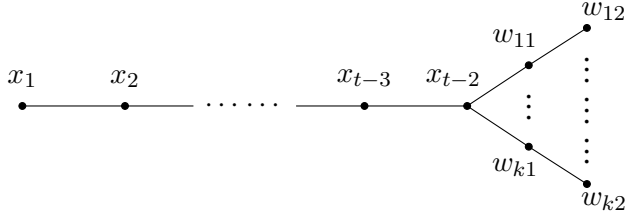
Fix a $t\ge 4$ and let $k$ be a positive integer. Consider the graph $T_t'$ with the vertices and the edges given as follows (see also \Cref{figure T'}):
\begin{align*}
    V(T_t')&=\{x_i,w_{j1},w_{j2}\mid i\in[t-2],j\in[k]\},\\
    E(T_t')&=\{\{x_i,x_{i+1}\},\{x_{t-2},w_{j1}\},\{w_{j1},w_{j2}\}\mid i\in[t-3],j\in [k]\}.
\end{align*}
From the structure of $T_t'$ one can easily see that $\left(I_t(T_t'):\prod_{i=1}^{t-2}x_i\right)=\l w_{j1}w_{j2}\mid j\in [k]\r$ and hence, by \Cref{reg sum}, $\reg\left(R/\left(I_t(T_t'):\prod_{i=1}^{t-2}x_i\right)\right)=k$. Thus, by \cite[Lemma 4.1]{CHHKTT}, $\reg(R/I_t(T_t'))\ge k$. On the other hand, $\nu_t(T_t')=1$.
\end{proof}

\begin{remark}\normalfont
    The above class of $t$-path ideals provides a counterexample to the recent conjecture by Hang and Vu \cite[Conjecture 4.9]{HangVu2024}.
\end{remark}

 Recall that, a graded ideal $I$ in $R$ is said to be {\it componentwise linear} if for each $j\ge 1$ the subideal $I_{\l j\r}$ generated by all degree $j$ elements in $I$ has regularity $\reg(I_{\l j\r})=j$. For a square-free monomial ideal $I\subseteq R$, the following implications hold due to 
\cite[Theorem 2.4]{MKA}, \cite[Chapter 8]{HHBook}, and Terai's formula:
\[
I^{\vee}\text{ is vertex splittable } \implies I^{\vee}\text{ is componentwise linear }\implies \pd(R/I)=\b(I).
\]
In \Cref{main theorem1}, we have shown that $I_3(T)^{\vee}$ is vertex splittable, which implies $I_3(T)^{\vee}$ is componentwise linear, and thus, we have $\pd(R/I_3(T))=\b(I_3(T))$ for any tree $T$. In fact, the above mentioned equality also follows from \Cref{pdmain}. On the other hand, it is well-known that the dual of the edge ideal (i.e., the dual of the $2$-path ideal) of a chordal graph is vertex splittable. In virtue of these results and trees being chordal graphs, it is natural to ask the following question:

\begin{question}\label{q1}
    Let $G$ be a chordal graph. Then, is it true that $I_{t}(G)^{\vee}$ is vertex splittable (or componentwise linear) for all $t\geq 3$? If not, then does there exist a chordal graph $G$ such that $I_t(G)^{\vee}$ is not vertex splittable (or componentwise linear), but $\pd(R/I_t(G))=\b(I_t(G))$, where $t\ge 3$?
\end{question}

We give an answer to the above question in the following theorem.

\begin{theorem}\label{Example1}
For each $t\ge 3$, there exists a connected chordal graph $G_t$ such that $I_t(G_t)^{\vee}$ is not componentwise linear, but $\pd(R/I_t(G_t))=\b(I_t(G_t))$.
\end{theorem}
\begin{proof}
    \begin{figure}[!ht]
\centering
\begin{tikzpicture}
[scale=.55]
\draw [fill] (0,2) circle [radius=0.1];
\draw [fill] (2,2) circle [radius=0.1];
\draw [fill] (5,2) circle [radius=0.1];
\draw [fill] (7,2) circle [radius=0.1];
\draw [fill] (8,3) circle [radius=0.1];
\draw [fill] (8,1)
circle [radius=0.1];
\draw [fill] (9,2) circle [radius=0.1];
\draw [fill] (11,2)
circle [radius=0.1];
\draw [fill] (14,2) circle [radius=0.1];
\draw [fill] (16,2) circle [radius=0.1];
\node at (0,2.5) {$x_1$};
\node at (2,2.5) {$x_2$};
\node at (3.5,2) {$\cdots$};
\node at (5,2.5) {$x_{t-2}$};
\node at (6.7,1.4) {$x_{t-1}$};
\node at (8,3.5) {$a$};
\node at (8,0.5) {$b$};
\node at (9,2.5) {$x_t$};
\node at (11,2.5) {$x_{t+1}$};
\node at (12.5,2) {$\cdots$};
\node at (14,2.5) {$x_{2t-3}$};
\node at (16,2.5) {$x_{2t-2}$};
\draw (0,2)--(2,2)--(3,2);
\draw (4,2)--(5,2)--(7,2)--(8,3)--(8,1)--(9,2)--(11,2)--(12,2);
\draw (13,2)--(14,2)--(16,2);
\draw (7,2)--(8,1);
\draw (8,3)--(9,2);
\end{tikzpicture}\caption{The graph $G_t$.}\label{figure 122466}
\end{figure}

Let $t\ge 3$ be an integer. Following \cite[Fig. 5]{ADGRS} consider the graph $G_t$ drawn as in \Cref{figure 122466}. From the structure of $G_t$ it is easy to see that for $W\subseteq V(G_t)$ with $|W|=t$ the induced subgraph $G_t[W]$ forms a $t$-path if and only if $G_t[W]$ is connected. Thus, $I_t(G_t)=J_t(G_t)$, where $J_t(G_t)$ is the $t$-connected ideal of $G_t$. Note that $J_t(G_t)$ is the Stanley-Reisner ideal of the $(t-1)$-independence complex $\mathrm{Ind}_{t-1}(G_t)$. Hence, by \cite[Proposition 4.3]{ADGRS}, $I_t(G_t)^{\vee}$ is not a componentwise linear ideal.

  We now proceed to show that $\pd(R/I_t(G_t))=\b(I_t(G_t))$. First note that if $\F$ forms a $t$-path in $G_t$, then either $a\in\F$ or $b\in\F$. Now consider the graph $G_t'$ with $V(G_t')=V(G_t)\setminus\{a\}$ and $E(G_t')=E(G_t\setminus a)\cup\{\{x_{t-1},x_t\}\}$. 
  
   \begin{figure}[!ht]
\centering
\begin{tikzpicture}
[scale=.55]
\draw [fill] (0,2) circle [radius=0.1];
\draw [fill] (2,2) circle [radius=0.1];
\draw [fill] (5,2) circle [radius=0.1];
\draw [fill] (7,2) circle [radius=0.1];
\draw [fill] (8,1)
circle [radius=0.1];
\draw [fill] (9,2) circle [radius=0.1];
\draw [fill] (11,2)
circle [radius=0.1];
\draw [fill] (14,2) circle [radius=0.1];
\draw [fill] (16,2) circle [radius=0.1];
\node at (0,2.5) {$x_1$};
\node at (2,2.5) {$x_2$};
\node at (3.5,2) {$\cdots$};
\node at (5,2.5) {$x_{t-2}$};
\node at (6.9,2.5) {$x_{t-1}$};
\node at (8,0.5) {$b$};
\node at (9,2.5) {$x_t$};
\node at (11,2.5) {$x_{t+1}$};
\node at (12.5,2) {$\cdots$};
\node at (14,2.5) {$x_{2t-3}$};
\node at (16,2.5) {$x_{2t-2}$};
\draw (0,2)--(2,2)--(3,2);
\draw (4,2)--(5,2)--(7,2)--(8,1)--(9,2)--(11,2)--(12,2);
\draw (13,2)--(14,2)--(16,2);
\draw (9,2)--(7,2)--(8,1);
\end{tikzpicture}\caption{The graph $G_t'$.}\label{figure G_t'}
\end{figure}
  
  Observe that for any $\F\subseteq V(G_t)$ with $a\notin\F$, $\F\cup\{a\}$ forms a $t$-path in $G_t$ if and only if $\F$ forms a $(t-1)$-path in $G_t'$. Also, if $\F$ forms a $t$-path in $G_t$ such that $a\notin\F$, then $(\F\setminus \{b\})\cup\{a\}$ forms a $t$-path in $G_t$, and hence, $\F\setminus\{b\}$ forms a $(t-1)$-path in $G_t'$. Therefore, $(I_t(G_t):a)=I_{t-1}(G_t')$. Moreover, $\l I_t(G_t),a\r=\l I_t(G_t\setminus a),a\r$, where $G_t\setminus a$ is a path graph of length $2t-2$. Hence, by \Cref{regularity lemma},
\begin{align}\label{eq1212}
\pd(R/I_t(G_t))\le\max\{\pd(R/I_{t-1}(G_t')),\pd(R/\l I_t(G_t\setminus a),a\r)\}.    
\end{align}
Our aim now is to show that $(I_{t-1}(G_t'):b)=I_{t-2}(G_t'\setminus\{x_1,b,x_{2t-2}\})$. Again, it is easy to see that, for any $\F\subseteq V(G_t')$ with $b\notin\F$, $\F\cup\{b\}$ forms a $(t-1)$-path of $G_t'$ if and only if $\F$ forms a $(t-2)$-path of $G_t'\setminus\{x_1,b,x_{2t-2}\}$. Now consider an $\F\subseteq V(G_t')$ with $b\notin \F$ such that $\F$ forms a $(t-1)$-path of $G_t'$. If $x_1\in \F$, then $(\F\setminus\{x_1\})\cup\{b\}$ forms a $(t-1)$-path of $G_t'$. Similarly,  if $x_{2t-2}\in\F$, then $(\F\setminus\{x_{2t-2}\})\cup\{b\}$ forms a $(t-1)$-path of $G_t'$. Now, for each $2\le i\le t-1$ if such an $\F$ forms the $(t-1)$-path of $G_t'$ which start from $x_i$, then $(\F\setminus\{x_{t+i-2}\})\cup\{b\}$ forms a $(t-1)$-path of $G_t'$. Similarly, for each $t\le i\le 2t-3$ if such an $\F$ forms the $(t-1)$-path of $G_t'$ which start from $x_i$, then $(\F\setminus\{x_{i-t+2}\})\cup\{b\}$ forms a $(t-1)$-path of $G_t'$. Thus, $(I_{t-1}(G_t'):b)=I_{t-2}(G_t'\setminus\{x_1,b,x_{2t-2}\})$. Hence, again by \Cref{regularity lemma},
\begin{align}\label{eq1213}
\pd(R/I_{t-1}(G_t'))\le\max\{\pd(R/I_{t-2}(G_t'\setminus\{x_1,b,x_{2t-2}\})),\pd(R/\l I_{t-1}(G_t'\setminus b),b\r)\}.
\end{align}
Now $G_t'\setminus\{x_1,b,x_{2t-2}\}$ is a path graph of length $2t-5$, and thus, $I_{t-2}(G_t'\setminus\{x_1,b,x_{2t-2}\})=J_{t-2}(G_t'\setminus\{x_1,b,x_{2t-2}\})$. Hence, by \cite[Thoerem 3.12]{ADGRS}, $I_{t-2}(G_t'\setminus\{x_1,b,x_{2t-2}\})$ is a vertex splittable ideal. Therefore, $\pd(R/I_{t-2}(G_t'\setminus\{x_1,b,x_{2t-2}\}))=\b(I_{t-2}(G_t'\setminus\{x_1,b,x_{2t-2}\}))=2$. Also, $G_t'\setminus b$ is a path graph of length of $2t-3$ and hence, $\pd(R/\l I_{t-1}(G_t'\setminus b),b\r)=1+\pd(R/I_{t-1}(G_t'\setminus b))=3$. Hence, by \Cref{eq1213}, $\pd(R/I_{t-1}(G_t'))\le 3$. Now $G_t\setminus a$ is a path graph of length $2t-2$ and hence, $\pd(R/\l I_t(G_t\setminus a),a\r)=1+\pd(R/ I_t(G_t\setminus a))=1+\b(I_t(G_t\setminus a))=3$. Therefore, using \Cref{eq1212} we get $\pd(R/I_t(G_t))=3$. It is easy to see that $\l x_2,a,x_{t+1}\r$ is a minimal prime ideal containing $I_t(G)$. Hence, $\pd(R/I_t(G_t))\le\b(I_t(G_t))$. Finally, using \Cref{lempdbght} we conclude that $\pd(R/I_t(G_t))=\b(I_t(G_t))$.
\end{proof}

In view of \Cref{Example1}, we see that \Cref{main theorem1} can not be generalized for chordal graphs. Now, in case of trees, it is natural to think about the generalization of \Cref{main theorem1} for all $t\geq 4$. That is, one may ask the following:

\begin{question}\label{qvstree}
   For any tree $T$, is $I_{t}(T)^{\vee}$ vertex splittable for all $t\ge 4$? If not, does the equality $\pd(R/I_t(T))=\b(I_t(T))$ hold for any tree $T$ and any $t\geq 4$?
\end{question}

The answer to the above question is negative even for the caterpillar graphs, which follows from the theorem given below.

\begin{theorem}\label{caterpillar not vs}
    For a caterpillar graph $G$, the difference between $\pd(R/I_t(G))$ and $\b(I_t(G))$ can be arbitrarily large for each $t\geq 4$.
\end{theorem}
\begin{proof}
     \begin{figure}[!ht]
\centering
\begin{tikzpicture}
[scale=.55]

\draw [fill] (1,2) circle [radius=0.1];
\draw [fill] (3,2) circle [radius=0.1];
\draw [fill] (7,2) circle [radius=0.1];
\draw [fill] (9,2) circle [radius=0.1];
\draw [fill] (0,0) circle [radius=0.1];
\draw [fill] (2,0) circle [radius=0.1];
\draw [fill] (8,0) circle [radius=0.1];
\draw [fill] (10,0) circle [radius=0.1];
\node at (1,2.5) {$x_1$};
\node at (3,2.5) {$x_2$};
\node at (6.8,2.5) {$x_{t-3}$};
\node at (9.2,2.5) {$x_{t-2}$};
\node at (4,2) {$\cdots$};
\node at (5,2) {$\cdots$};
\node at (6,2) {$\cdots$};
\node at (1,0) {$\cdots$};
\node at (9,0) {$\cdots$};
\node at (-0.2,-0.5) {$y_1$};
\node at (2.2,-0.5) {$y_n$};
\node at (7.8,-0.5) {$z_1$};
\node at (10.2,-0.5) {$z_m$};
\draw (0,0)--(1,2)--(3,2);
\draw (2,0)--(1,2);
\draw (7,2)--(9,2)--(8,0);
\draw (10,0)--(9,2);
\end{tikzpicture}\caption{The graph $T_t$.}\label{figure T}
\end{figure}
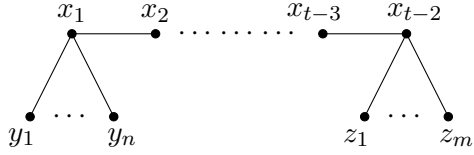
Fix a $t\ge 4$ and let $n,m$ be two positive integers. Consider the caterpillar graph $T_t$ with the vertex set and edge set given as follows (see also \Cref{figure T}):
\begin{align*}
    V(T_t)&=\{x_i,y_j,z_k\mid i\in[t-2],j\in[n],k\in [m]\},\\
    E(T_t)&=\{\{x_1,y_j\},\{x_{t-2},z_k\},\{x_i,x_{i+1}\}\mid j\in[n],k\in[m],i\in [t-3]\}.
\end{align*}
Analyzing the structure of $T_t$ one can easily see that 
\[
I(T_t)^{\vee}=\left\l x_1,\ldots,x_{t-2},\prod_{j=1}^ny_j,\prod_{k=1}^mz_k \right\r.
\]
Therefore, by \Cref{reg sum}, $\reg(I(T_t)^{\vee})=n+m-1$, and hence, by Terai's formula $\pd(R/I_t(T))=n+m-1$. However, $\b(I_t(T))=\max\{n,m\}$.
\end{proof}

\subsection{Concluding remarks:}\label{conclusion} 
As discussed earlier, the $t$-path ideals of graphs are the same as the $t$-connected ideals of graphs if and only if $t=2,3$. So far, we have seen that most of our main results proved earlier for $3$-path ideals of chordal graphs can not be generalized to higher $t$-path ideals, even for the class of trees. Now, it is natural to ask analogous questions of \Cref{qreg}, \ref{q1}, and \ref{qvstree} in the context of $t$-connected ideals. Regarding this, we have the following observations:
\begin{enumerate}

    \item[$\bullet$] In \Cref{reg tree}, we have shown that \Cref{result1} and \Cref{caterpillar reg} can not be generalized for arbitrary tree when $t\geq 4$, i.e., the answer to \Cref{qreg} is negative for $t$-path ideals. However, the $t$-connected ideals of the trees considered in \Cref{reg tree} has a linear resolution by \cite[Theorem 5.1]{AJM2024} as $\nu(\H)=1$, where $\H$ is the hypergraph corresponding to the $t$-connected ideal $J_t(T_{t}')$. In other words, $\reg(R/J_t(T_{t}'))=t-1=(t-1)\nu(\H)$. 

    \item[$\bullet$] The answer to \Cref{q1} was provided in \Cref{Example1}. As mentioned earlier, for the graphs considered in \Cref{Example1}, $t$-path ideals and $t$-connected ideals coincide. Thus, we have a similar answer to the analogue of \Cref{q1} for $t$-connected ideals as compared to the $t$-path ideals.

    \item[$\bullet$] In \Cref{caterpillar not vs}, we give examples of caterpillar graphs to show the answer of \Cref{qvstree} is negative. But, for any caterpillar graph $G$, it follows from \cite[Theorem 3.12]{ADGRS} that $J_{t}(G)^{\vee}$ is vertex splittable for all $t\geq 2$. Till now, we have no evidence to counter the analogue of \Cref{qvstree} for $t$-connected ideals.
\end{enumerate}

In view of the above discussions, it is reasonable to believe that the notion of $t$-connected ideals might be a more suitable generalization of edge ideals than $t$-path ideals, while dealing with homological invariants like regularity and projective dimension. Due to all of these facts and our computational evidence, we propose the following conjectures.

\begin{conjecture}
    Let $T$ be a tree and let $J_t(T)$ denote the $t$-connected ideal of $T$. Then $\reg(R/J_t(T))=(t-1)\nu(\H)$, where $\H$ is the hypergraph corresponding to the  $t$-connected ideal $J_t(T)$.
\end{conjecture}

\begin{conjecture}
   Let $T$ be a tree and let $J_t(T)$ denote the $t$-connected ideal of $T$. Then $\pd(R/J_t(T))=\b(J_t(T))$.
\end{conjecture}

A homogeneous ideal $I$ generated in degree $d$ is said to have a linear resolution if $\reg(R/I)=d-1$. Fr\"oberg \cite{Froberg1990} characterized all quadratic square-free monomial ideals with linear resolution. However, for arbitrary square-free monomial ideals, such a characterization is wide open. For the $3$-path ideals of chordal graphs, it follows from \Cref{result1} that $I_3(G)$ has a linear resolution if and only if $\nu_3(G)=1$. Also, when $G$ is co-chordal, it has been proved in \cite[Theorem 3.12]{DRSV23} that $I_3(G)$ has a linear resolution. Note that the $5$-cycle $C_5$ is a graph which is neither chordal nor co-chordal, but $I_3(C_5)$ has a linear resolution \cite[Corollary 5.5]{AlilooeeFaridi2015}. In fact, from the recent results of Hang and Vu \cite{HangVu2024}, we get some unicyclic graphs $G$ for which $I_3(G)$ has a linear resolution. Thus, it is worthwhile to explore the following question:

\begin{question}
    Characterize all graphs $G$ for which $I_3(G)$ has a linear resolution.
\end{question}

\noindent
{\bf Acknowledgements.} The first and the second authors are supported by Postdoctoral Fellowships at Chennai Mathematical Institute. The third author would like to thank the National Board for Higher Mathematics (India) for the financial support through the NBHM Postdoctoral Fellowship. All the authors are partially supported by a grant from the Infosys Foundation.  The authors acknowledge the use of the computer algebra system Macaulay2 \cite{M2} and the online platform SageMath for testing their computations.

\bibliographystyle{abbrv}
\bibliography{ref}

\end{document}